\documentclass[reqno,11pt]{amsart} 
\usepackage{amsmath}
\usepackage{amssymb}
\usepackage{amsthm}

\newcommand\NoBlackBoxes{\global\overfullrule0pt}
\NoBlackBoxes
\theoremstyle{plain}

\begin{document}

\title{Asymptotic behavior of R\'enyi entropy \\ 
in the central limit theorem \\
{\rm \tiny (extended version)}
}

\author{Sergey G. Bobkov$^{1}$}
\thanks{1) 
School of Mathematics, University of Minnesota, Minneapolis, MN 55455 USA.
Research was partially supported by NSF grant DMS-1612961.}

\author{Arnaud Marsiglietti$^{2}$}
\thanks{2) 
Center for the Mathematics of Information, California Institute of Technology, 
Pasadena, CA 91125 USA. 
\newline \hskip10mm
Supported by the Walter S. Baer and Jeri Weiss CMI 
Postdoctoral Fellowship}


\begin{abstract}
We explore an asymptotic behavior of R\'enyi entropy along convolutions in 
the central limit theorem with respect to the increasing number of i.i.d. 
summands. In particular, the problem of monotonicity is addressed under
suitable moment hypotheses.
\end{abstract}

\maketitle
\markboth{Sergey G. Bobkov and Arnaud Marsiglietti}{
Expansions for R\'enyi entropy in the CLT}

\def\theequation{\thesection.\arabic{equation}}
\def\E{{\mathbb E}}
\def\R{{\mathbb R}}
\def\C{{\mathbb C}}
\def\P{{\mathbb P}}
\def\H{{\rm H}}
\def\Im{{\rm Im}}
\def\Tr{{\rm Tr}}

\def\k{{\kappa}}
\def\M{{\cal M}}
\def\Var{{\rm Var}}
\def\Ent{{\rm Ent}}
\def\O{{\rm Osc}_\mu}

\def\ep{\varepsilon}
\def\phi{\varphi}
\def\vp{\varphi}
\def\F{{\cal F}}
\def\L{{\cal L}}

\def\be{\begin{equation}}
\def\en{\end{equation}}
\def\bee{\begin{eqnarray*}}
\def\ene{\end{eqnarray*}}

\vskip5mm
\section{{\bf Introduction}}
\setcounter{equation}{0}

\vskip2mm
\noindent
Given a (continuous) random variable $X$ with density $p$, the associated 
R\'enyi entropy and R\'enyi entropy power of index $r$ ($1 < r < \infty$)
are defined by
$$
h_r(X) \, = \, -\frac{1}{r-1}\,\log \int_{-\infty}^\infty p(x)^r\,dx, \qquad
N_r(X) \, = \, e^{2h_r(X)} \, = \,
\bigg(\int_{-\infty}^\infty p(x)^r\,dx\bigg)^{-\frac{2}{r-1}}.
$$
Being translation invariant and homogeneous of order 2, the functional 
$N_r$ is similar to the variance and is often interpreted as measure 
of uncertainty hidden in the distribution of $X$. Another representation
$$
N_r(X)^{-\frac{1}{2}} \, = \, \big(\E\, p(X)^{r-1}\big)^{\frac{1}{r-1}}
$$
shows that $N_r$ is non-increasing in $r$,
so that $0 \leq N_\infty \leq N_r \leq N_1 \leq \infty$.
Here, for the extreme indexes, the R\'enyi entropy power is defined by the
monotonicity,
$$
N_\infty(X) = \lim_{r \uparrow \infty} N_r(X) = \|p\|_\infty^{-2}, \qquad
N_1(X) = \lim_{r \downarrow 1} N_r(X) = e^{2h_1(X)},
$$
where $\|p\|_\infty$ is the essential supremum of $p(x)$. In the case $r=1$,
we arrive at the Shannon differential entropy 
$h_1(X) = h(X) = -\int p(x)\log p(x)\, dx$ with entropy power 
$N_1 = N = e^{2h}$
(provided that $N_r(X) > 0$ for some $r>1$).

Much of the analysis about the Shannon and R\'enyi entropies is focused
on the behavior of these functionals on convolutions, i.e.,
for sums $S_n = X_1+\dots+X_n$ of independent random variables
(including a multidimensional setting). First, let us recall
a fundamental entropy power inequality, which may be written in terms
of the normalized sums $Z_n = S_n/\sqrt{n}$~as
\be
N(Z_n) \, \geq \, \frac{1}{n}\,\sum_{k=1}^n N(X_k).
\en
There are also some extensions of this relation to the R\'enyi case 
(cf. [D-C-T], [C-T], [B-C], [B-M]). 

When $X_k$'s are independent and identically distributed (i.i.d.), with 
mean zero and variance one, the central limit theorem (CLT) asserts that
$Z_n \Rightarrow Z$ with weak convergence in distribution to the Gaussian 
limit $Z \sim N(0,1)$.
In this case, the right-hand side of (1.1) is constant, while the sequence 
on the left is monotone, as was shown by Artstein, Ball, Barthe and Naor 
[A-B-B-N], cf. also [M-B] (the inequality (1.1) itself ensures that $N(Z_n)$ 
are non-decreasing along only the power values $n = 2^l$). Moreover, by another 
important result due to Barron [B], we have the entropic CLT: $N(Z_n)$ 
are convergent to the entropy power $N(Z)$, as long as $N(Z_{n_0}) > 0$ 
for some $n_0$. 

These results give rise to a number of natural questions about an asymptotic 
behavior of the R\'enyi entropy powers $N_r(Z_n)$. In particular, when 
do they converge to $N_r(Z)$, and if so, what is the rate of convergence?
Is the monotonicity still true? As we will see, such questions may be
studied, at least partially, under suitable moment assumptions.

Let us state a few observations in these directions, assuming throughout 
that $X,X_1,X_2,\dots$ are i.i.d. random variables with $\E X = 0$ and
$\Var(X) = 1$. Put $\beta_s = \E\,|X|^s$ for real $s \geq 2$. In order to 
describe necessary and sufficient conditions for the convergence of the R\'enyi entropies in the CLT, we also introduce the common characteristic function
$$
f(t) = \E\,e^{itX} \qquad (t \in \R).
$$

\vskip5mm
{\bf Theorem 1.1.} {\sl Given $1 < r \leq \infty$, we have the convergence 
$N_r(Z_n) \rightarrow N_r(Z)$ or equivalently 
$h_r(Z_n) \rightarrow h_r(Z)$ as $n \rightarrow \infty$, if and only if 
\be
\int_{-\infty}^\infty |f(t)|^\nu\,dt < \infty \quad {\sl for \ some} \ \ \nu \geq 1.
\en
Equivalently, this holds if and only if $Z_n$ have bounded densities for all 
$($some$)$ $n$ large enough.
}

\vskip5mm
This characterization coincides with the one for the uniform local limit 
theorem due to Gnedenko, cf. [G-K]. Since (1.2) is equivalent to the property 
that $Z_n$ have bounded and hence bounded $C^k$-smooth densities for any 
fixed $k$ and all $n$ large enough, it is often referred to as the
smoothing condition. In general, (1.2) is stronger than what is needed
in the entropic case $r=1$. In this connection, let us note that there is 
still no explicit description such as (1.2) for the validity of the entropic 
CLT in terms of the characteristic function $f(t)$.

Once (1.2) is fulfilled, one may ask about the rate of convergence in 
Theorem 1.1, which may be guaranteed assuming that the absolute moment $\beta_s$ 
is finite for some $s>2$. Moreover, in this case one may obtain asymptotic
expansions for $N_r(Z_n)$ in powers of $1/n$ similarly to the entropic 
expansions derived in [B-C-G2]. They involve the moments of $X$ up to order 
$m = [s]$, or equivalently -- the cumulants
$$
\gamma_k = i^{-k}\,(\log f)^{(k)}(0), \qquad k = 1,\dots,m.
$$
In the Gaussian case $X \sim N(0,1)$, all cumulants are vanishing,
starting with $k = 2$. In the general case, they indicate how close 
a given distribution to the normal. In the asymptotic behavior of R\'enyi's
entropies, it turns out that a special role is played by the quantity
$$
b = b(r) = -\frac{1}{r}\ \bigg[
\frac{2-r}{12}\, \gamma_3^2 + \frac{r-1}{8}\, \gamma_4\bigg].
$$
Here, $\gamma_3 = \E X^3$ and $\gamma_4 = \E X^4 - 3$, while for the extreme
indexes, one may just put
$$
b(1) = \lim_{r \rightarrow 1} b(r) = -\frac{1}{12}\,\gamma_3^2, \qquad
b(\infty) = \lim_{r \rightarrow \infty} b(r) =
\frac{1}{12}\, \gamma_3^2 - \frac{1}{8}\, \gamma_4.
$$
This can be seen from the following assertion.

\vskip5mm
{\bf Theorem 1.2.} {\sl Suppose that the smoothing condition $(1.2)$ is
fulfilled. If $\beta_s$ is finite for $2 \leq s < 4$, then 
for any $1 < r < \infty$,
\be
h_r(Z_n) = h_r(Z) + o(n^{-\frac{s-2}{2}}), \qquad
N_r(Z_n) = N_r(Z) + o(n^{-\frac{s-2}{2}}).
\en
Moreover, in case $4 \leq s < 6$,
\begin{eqnarray}
h_r(Z_n) 
 & = &
h_r(Z) + b\, n^{-1} + o(n^{-\frac{s-2}{2}}),\\
N_r(Z_n) 
 & = &
N_r(Z)\,\big(1 + 2 b\, n^{-1}\big) + o(n^{-\frac{s-2}{2}}) \nonumber .
\end{eqnarray}
}

\vskip2mm
This assertion remains valid in the entropic case $r=1$ as well
(with a slight logarithmic improvement in the remainder $o$-term, cf. [B-C-G2]).
In case $s=6$, the remainder term may be improved to $O(n^{-2})$, and in fact,
one may add quadratic terms to get an expansion
\be
h_r(Z_n)  \, = \,
h_r(Z) + b\, n^{-1} + b_2 n^{-2} + o(n^{-2})
\en
with some functional $b_2 = b_2(r)$ depending also on $\gamma_5$ and $\gamma_6$. 
Regardless of its value, 
one may therefore conclude about an eventual monotonicity of $N_r(Z_n)$
based on the sign of $b$. Moreover, the above expansions continue to hold 
for $r = \infty$, so that this case may be included as well.

\vskip5mm
{\bf Theorem 1.3.} {\sl Suppose that the smoothing condition $(1.2)$ is
fulfilled, and let $\beta_6$ be finite. Given $1 < r \leq \infty$, there 
exists $n_0 \geq 1$ such that the sequence $N_r(Z_n)$ is increasing 
for $n \geq n_0$, whenever $b(r) < 0$, that is, if
$$
\frac{2-r}{3}\, \gamma_3^2 + \frac{r-1}{2}\, \gamma_4 > 0 \quad
(1 < r < \infty), \qquad
\gamma_4 > \frac{2}{3}\, \gamma_3^2 \quad (r = \infty).
$$
This sequence is decreasing for $n \geq n_0$, if $b(r) > 0$.
}

\vskip5mm
In particular, under the last condition $\gamma_4 > \frac{2}{3}\, \gamma_3^2$, 
the sequence $N_r(Z_n)$ is eventually increasing for any fixed $r \geq 1$.
For example, this holds for $X = \frac{\xi - \alpha}{\sqrt{\alpha}}$,
where the random variable $\xi$ has a Gamma distribution with $\alpha$ degrees
of freedom (in which case $\gamma_3 = 2/\sqrt{\alpha}$ and 
$\gamma_4 = 6/\alpha$). 

On the other hand, if $X$ is uniformly distributed
in the interval $(-\sqrt{3},\sqrt{3})$, then $\gamma_3 = 0$, $\gamma_4 = -6/5$, 
so $N_r(Z_n)$ is eventually decreasing for any $r > 1$, although 
the opposite property takes place for $r=1$.

The paper is organized as follows. We start with the proof of Theorem 1.1
(Section 2), and then collect together basic results on Edgeworth expansions for 
densities $p_n$ of $Z_n$ (Section~3). They are used in Sections 4-5 to construct 
a formal asymptotic expansion for $L^r$-norms of $p_n$ in powers of $1/n$ up to 
order $[\frac{m-2}{2}]$ with remainder term as in (1.3)-(1.4). One particular 
case, where the first moments of $X$ agree with those of $Z \sim N(0,1)$, is 
discussed separately in Section~6, while the range $4 \leq s \leq 8$ in such 
expansion is treated in Section 7. The transition to the R\'enyi entropy 
is performed in Section 8, where Theorem 1.2 is proved. Some comparison with 
the entropic CLT is given in Section 9, with remarks leading to Theorem 1.3
for finite $r$. Finally, the index $r = \infty$ is treated separately in 
Section 10. We thus follow the next plan:

\vskip5mm
1. Introduction

2. Proof of Theorem 1.1

3. Limit theorems about Edgeworth expansions

4. Approximation for $L^r$-norm of densities $p_n$

5. Truncated $L^r$-norm of approximating densities $\varphi_m$

6. The case where the first cumulants are vanishing

7. Moments of order $4 \leq s \leq 8$

8. Expansions for R\'enyi entropies

9. Comparison with the entropic CLT. Monotonicity

10. Maximum of density (the case $r = \infty$)

\vskip5mm
\section{{\bf Proof of Theorem 1.1}}
\setcounter{equation}{0}

\vskip2mm
\noindent
From now on, let $X,X_1,X_2,\dots$ be i.i.d. random variables with $\E X = 0$ 
and $\Var(X) = 1$, for which we define the normalized sums
$$
Z_n = \frac{X_1 + \dots + X_n}{\sqrt{n}}, \qquad n = 1,2,\dots
$$

First, let us recall Gnedenko's uniform local limit theorem. Assuming the 
smoothing condition (1.2), it asserts that, for all $n$ large enough, the random
variables $Z_n$ have bounded densities $p_n$, and moreover, in that case as 
$n \rightarrow \infty$,
\be
\sup_x |p_n(x) - \varphi(x)| \rightarrow 0.
\en
Here, as usual,
$$
\varphi(x) = \frac{1}{\sqrt{2\pi}}\,e^{-x^2/2} \qquad (x \in \R)
$$ 
denotes the density of the standard normal random variable $Z$. Clearly,
the property (2.1) is also necessary for the uniform boundedness of $p_n$'s.

Let us explain the equivalence of the two conditions -- in terms of the 
characteristic function as in (1.2), and in terms of densities (via the existence 
of a bounded density). Since $|f(t)| \leq 1$ for all $t$, the property (1.2) 
is getting weaker for growing $\nu$, so it is sufficient to consider integer 
values of $\nu$. Since $Z_n$ has characteristic function
$$
f_n(t) = \E\,e^{itZ_n} = f(t/\sqrt{n})^n,
$$
(1.2) implies that $Z_n$ has a bounded, continuous density $p_n$ for $n = \nu$, 
by the Fourier inversion formula. Hence the same is true for all $n \geq \nu$, 
by the convolution character of the distributions of $Z_n$. 
Conversely, suppose that $Z_n$ has a bounded density $p_n$ for $n = n_0$. 
This implies that $p_n \in L^r(\R)$ for any $r \geq 1$, with norm
$$
\|p_n\|_r = \bigg(\int_{-\infty}^\infty p_n(x)^r\,dx\bigg)^{1/r},
$$
and in particular
$p_n \in L^2(\R)$. By Plancherel's theorem, the characteristic function $f_n$ 
is also in $L^2(\R)$. But this means that (1.2) is fulfilled 
with $\nu = 2n_0$.

Also note that, under the condition (1.2), we have $f_\nu(t) \rightarrow 0$
as $t \rightarrow \infty$ (the Riemann-Lebesgue lemma), and thus 
$f(t) \rightarrow 0$. Hence, (1.2) represents a sharpening of 
the Cram\'er condition $\limsup_{t \rightarrow \infty} |f(t)| < 1$, 
which is used to establish a number of asymptotic results related to 
the CLT. In particular, using the Fourier inversion formula, 
one can easily obtain (2.1) and actually a sharper statement such as
\be
\sup_x \  (1+x^2)\,|p_n(x) - \varphi(x)| \rightarrow 0 \qquad 
(n \rightarrow \infty).
\en

\vskip2mm
{\bf Proof of Theorem 1.1.} First, let $r = \infty$.
As explained, the smoothing condition (1.2) implies the uniform local
limit theorem (2.1). In turn, the latter yields 
$\|p_n\|_\infty \rightarrow \|\varphi\|_\infty$, that is,
$N_\infty(Z_n) \rightarrow N_\infty(Z)$ as $n \rightarrow \infty$.
Conversely, this convergence ensures that $N_\infty(Z_n) > 0$ for all $n$ 
large enough, that is, $\|p_n\|_\infty < \infty$. As was also emphasized, 
this implies (1.2).

Now, let $1 < r < \infty$. In one direction, if $N_r(Z_n) \rightarrow N_r(Z)$ 
as $n \rightarrow \infty$, then $N_r(Z_n) > 0$ for all $n$ large enough, say 
$n \geq n_0$. Equivalently, for such $n$, $Z_n$ have densities $p_n$ with 
$\|p_n\|_r < \infty$. If $r \geq 2$, then $\|p_n\|_2 \leq 1 + \|p_n\|_r < \infty$, 
so that $p_n$ and therefore $f_n$ are in $L^2(\R)$. This means that (1.2) is 
fulfilled for $\nu = 2n_0$. In the case $1 < r < 2$, one may apply 
the Hausdorff-Young inequality
$$
\|\hat u\|_{r'} \leq \|u\|_{r}, \qquad {\rm where} \quad
\hat u(t) = \int_{-\infty}^\infty e^{2\pi i tx}\,u(x)\,dx, \quad
r' = \frac{r}{r-1}.
$$
It implies that $\|f_n\|_{r'} \leq \sqrt{2\pi}\,\|p_n\|_{r} < \infty$, 
which means that (1.2) is fulfilled for $\nu = r' n_0$.

Thus, the smoothing condition (1.2) is indeed necessary. To argue in the other
direction, we apply the uniform local limit theorem: For all $n \geq n_0$ 
large enough, $Z_n$ have densities $p_n$, bounded by a constant $M$ and moreover, 
the relation (2.1) holds true, i.e.,
\be
\sup_x \big|p_n(x)^r - \varphi(x)^r\big| \leq  \ep_n \rightarrow 0 
\qquad (n \rightarrow \infty).
\en
For a given $\ep > 0$, applying the usual central limit theorem, one may pick up 
$T>0$ such that
$$
\P\{|Z_n| > T\} + \P\{|Z| > T\} < \ep, \qquad n \geq n_1 \geq n_0.
$$
Hence
$$
\int_{|x|>T} p_n(x)^r\,dx \, \leq \, M^{r-1} \int_{|x|>T} p_n(x)\,dx \, = \,
M^{r-1}\, \P\{|Z_n| > T\} \, < \, M^{r-1} \ep,
$$
and similarly for $\varphi(x)$. Hence
\be
\bigg|\int_{|x|>T} p_n(x)^r\,dx  - \int_{|x|>T} \varphi(x)^r\,dx\bigg| < 
M^{r-1} \ep.
\en
On the other hand, by (2.3),
\bee
\bigg|\int_{|x|\leq T} p_n(x)^r\,dx  - \int_{|x| \leq T} \varphi(x)^r\,dx\bigg| 
 & \leq &
\int_{|x| \leq T} |p_n(x)^r  - \varphi(x)^r|\,dx \\
 & \leq &
2T \ep_n \ \leq \ \ep,
\ene
where the last inequality holds true for all $n \geq n_2$ with some $n_2 \geq n_1$. 
Together with (2.4), we get
$$
\big|\,\|p_n\|_r^r - \|\varphi\|_r^r\big| < (M^{r-1} + 1)\, \ep, \qquad 
n \geq n_2.
$$
That is, $\|p_n\|_r^r \rightarrow \|\varphi\|_r^r$ as $n \rightarrow \infty$,
thus proving the theorem.
\qed

\vskip5mm
\section{{\bf Limit Theorems about Edgeworth Expansions}}
\setcounter{equation}{0}

\vskip2mm
\noindent
As is well-known, in case of the finite 3-rd absolute moment 
$\beta_3 = \E\,|X|^3$, and assuming the smoothness condition (1.2), the local 
limit theorem (2.1) or even the non-uniform variant (2.2) can be sharpened to
\be
\sup_x \  (1 + |x|^3)\,|p_n(x) - \varphi(x)| \, = \,
o\bigg(\frac{1}{\sqrt{n}}\bigg) \qquad (n \rightarrow \infty).
\en
Here, the rate cannot be improved in general. However, under higher order moment
assumptions, the limit normal density may slightly be modified, which leads 
to the sharpening of the right-hand side of (3.1). Namely, if 
$\beta_m = \E\,|X|^m$ is finite
for an integer $m \geq 2$, one may introduce the cumulants
$$
\gamma_k = i^{-k}\,(\log f)^{(k)}(0), \qquad k = 1,\dots,m.
$$
They represent certain polynomials in the moments $\alpha_i = \E X^i$ up 
to order $k$, namely,
$$
\gamma_k \, = \, k!\, \sum \ (-1)^{j - 1}\, (j-1)! \,\frac{1}{r_1! \dots r_k!} \ 
\Big(\frac{\alpha_1}{1!}\Big)^{r_1} \dots
\Big(\frac{\alpha_k}{k!}\Big)^{r_k},
$$
where $j = r_1 + \dots + r_k$ and where the summation is running over all 
tuples $(r_1,\dots,r_k)$ of non-negative integers such that 
$r_1 + 2 r_2 + \dots + k r_k = k$.

For example, with our moment assumptions $\E X = 0$, 
$\Var(X) = 1$, we have $\gamma_1 = 0$, $\gamma_2 = 1$,
$$
\gamma_3 = \alpha_3, \quad \gamma_4 = \alpha_4 - 3, \quad 
\gamma_5 = \alpha_5 - 10\, \alpha_3, \qquad 
\gamma_6 = \alpha_6 - 15\,\alpha_4 - 10\,\alpha_3^3 + 30.
$$ 

\vskip2mm
{\bf Definition 3.1.} An Edgeworth correction of the standard normal law
of order $m$ for the distribution of $Z_n$ is a finite signed measure $\nu_m$ 
with density
\be
\varphi_m(x) = \varphi(x) + \varphi(x)\, \sum_{k=1}^{m-2} Q_k(x)\,n^{-k/2},
\en
where
\be
Q_k(x) = \sum \frac{1}{r_1! \dots r_k!} \ 
\Big(\frac{\gamma_3}{3!}\Big)^{r_1} \dots
\Big(\frac{\gamma_{k+2}}{(k+2)!}\Big)^{r_k} \ H_{k+2j}(x).
\en
Here, the summation is running over all collections 
of non-negative integers $r_1,\dots,r_k$ such that 
$r_1 + 2r_2 + \dots + k r_k = k$, with notation $j = r_1 + \dots + r_k$.

\vskip5mm
As usual, $H_k$ denotes the Chebyshev-Hermite polynomial of degree $k$ 
with leading term $x^k$. The polynomial $Q_k$ in (3.2) has degree at most 
$3(m-2)$ in the variable $x$. Indeed, the index 
$$
k+2j = 3r_1 + 4r_2 + \dots + (k+2) r_k = k + 2(r_1 + \dots + r_k)
$$
is maximized for $k=m-2$ and for the collection 
$r_1 = m-2$, $r_2 = \dots = r_{m-2} = 0$. In this case, (3.3) contains the term
$$
\frac{1}{(m-2)!}\,\Big(\frac{\gamma_3}{3!}\Big)^{m-2}\, H_{3(m-2)}(x).
$$
of degree exactly $3(m-2)$ as long as $\gamma_3 \neq 0$.

The index $m$ for $\varphi_m$ indicates that the cumulants up to $\gamma_m$
participate in the construction. The sum in (3.2) may also be viewed
as a polynomial in $1/\sqrt{n}$ of degree at most $m-2$.

For example, $\varphi_2 = \varphi$, and there are no terms in the sum (3.2).
For $m=3,4,5,6$, in (3.3) we correspondingly have
\bee
Q_1(x) 
 & = &
\frac{\gamma_3}{3!}\, H_3(x), \\
Q_2(x) 
 & = &
\frac{\gamma_3^2}{2!\,3!^2}\, H_6(x) + \frac{\gamma_4}{4!}\, H_4(x), \\
Q_3(x) 
 & = &
\frac{\gamma_3^3}{3!^4}\, H_9(x) + \frac{\gamma_3 \gamma_4}{3!\,4!}\, H_7(x) + 
\frac{\gamma_5}{5!}\, H_5(x), \\
Q_4(x) 
 & = &
\frac{\gamma_3^4}{4!\,3!^4}\, H_{12}(x) + 
\frac{\gamma_3^2 \gamma_4}{2!\, 3!^2\,4!}\, H_{10}(x) + 
\frac{\gamma_3 \gamma_5}{3!\,5!}\, H_8(x) + 
\frac{\gamma_4^2}{2!\, 4!^2}\, H_8(x) + 
\frac{\gamma_6}{6!}\, H_6(x). 
\ene

Moreover, if the first $m-1$ moments of $X$ coincide with those of
$Z \sim N(0,1)$, then the first $m-1$ cumulants of $X$ are vanishing, and (3.2)
is simplified to
\be
\varphi_m(x) = \varphi(x)\bigg(1 + \frac{\gamma_m}{m!}\, H_m(x)\, 
n^{-\frac{m-2}{2}}\bigg),
\en
where necessarily $\gamma_m = \E X^m - \E Z^m$.

The following observation, generalizing and refining the non-uniform 
local limit theorems (2.2) and (3.1), is due to Petrov [P1], 
cf. also [P2], [B-RR]. From now on, we always assume that the smoothing
condition (1.2) is fulfilled.

\vskip5mm
{\bf Lemma 3.2.} {\sl If $\beta_m < \infty$ for an integer $m \geq 2$, then
as $n \rightarrow \infty$
\be
\sup_x \  (1 + |x|^m)\,|p_n(x) - \varphi_m(x)| \, = \,
o\big(n^{-\frac{m-2}{2}}\big).
\en
}

\vskip2mm
Without the polynomial weight $1 + |x|^m$, a similar result was earlier 
obtained by Gnedenko. However, in some applications the appearance of 
this weight turns out to be crucial.

If $m \geq 3$, one may also take $\varphi_{m-1}$ as 
an approximation of $p_n$, and then (3.5) together with Definition 3.1 imply that
\be
\sup_x \  (1 + |x|^m)\,|p_n(x) - \varphi_{m-1}(x)| \, = \,
O\big(n^{-\frac{m-2}{2}}\big).
\en

A further generalization was given in [B-C-G1] to employ the case
of fractional moments.

\vskip5mm
{\bf Lemma 3.3.} {\sl Let $\beta_s < \infty$ for some real $s \geq 2$, 
and $m = [s]$. Then uniformly over all $x$, as $n \rightarrow \infty$,
$$
(1 + |x|^s)\,(p_n(x) - \varphi_m(x)) \, = \,
o\big(n^{-\frac{s-2}{2}}\big) + (1 + |x|^{s-m})\,
\Big(O\big(n^{-\frac{m-1}{2}}\big) + o\big(n^{-(s-2)}\big)\Big).
$$
In particular, for some constant $\alpha>0$ depending on $s$,
\be
\sup_{|x| \leq n^\alpha} \  (1 + |x|^s)\,|p_n(x) - \varphi_m(x)| \, = \,
o\big(n^{-\frac{s-2}{2}}\big).
\en
}

Thus, (3.7) extends (3.5) when taking the supremum over relatively 
large interval.

There are also similar results about the distribution functions
$F_n(x) = \P\{Z_n \leq x\}$, which may be approximated by
\be
\Phi_m(x) \, = \, \nu_m((-\infty,x]) \, = \, \int_{-\infty}^x \varphi_m(y)\,dy 
\, = \, \Phi(x) - \varphi(x)\, \sum_{k=1}^{m-2} R_k(x)\,n^{-k/2},
\en
where
$$
R_k(x) = \sum \frac{1}{r_1! \dots r_k!} \ 
\Big(\frac{\gamma_3}{3!}\Big)^{r_1} \dots
\Big(\frac{\gamma_{k+2}}{(k+2)!}\Big)^{r_k} \ H_{k+2j-1}(x)
$$
with summation as in Definition 3.1. The next result is due to Osipov and
Petrov [O-P].

\vskip5mm
{\bf Lemma 3.4.} {\sl Suppose that $\beta_s < \infty$ for some real $s \geq 2$, 
and let $m = [s]$. Then, as $n \rightarrow \infty$,
$$
\sup_x \
(1 + |x|^s)\,|F_n(x) - \Phi_m(x)| \, = \, o\big(n^{-\frac{s-2}{2}}\big).
$$
In particular, when $s = m \geq 3$ is integer, we have
$$
\sup_x \
(1 + |x|^s)\,|F_n(x) - \Phi_{m-1}(x)| \, = \, O\big(n^{-\frac{s-2}{2}}\big).
$$
}

This statement holds under the weaker assumption in comparison
with (1.2): nothing should be required in case $2 \leq s < 3$, while for 
$s \geq 3$ the Cram\'er condition is sufficient.

\vskip5mm
{\bf Remark 3.5.}
Since the densities $p_n$ can properly be approximated by the functions 
$\varphi_m$, it makes sense to isolate the leading term in the sum (3.2), 
by rewriting the definition as
\be
\varphi_m(x) = \varphi(x) + 
\varphi(x) \frac{\gamma_{k+2}}{(k+2)!}\, H_{k+2}(x)\,n^{-k/2} + 
\varphi(x)\, \sum_{j=k+1}^{m-2} Q_j(x)\,n^{-j/2}
\en
for some unique $1 \leq k \leq m-2$. The value of $k$ is the maximal one 
in the interval $[1,m-2]$ such that $\gamma_3 = \dots = \gamma_{k+1} = 0$,
which means that the first moments of $X$ up to order $k+1$ coincide with those of
$Z \sim N(0,1)$. In this case, necessarily 
$$
\gamma_{k+2} = \E X^{k+2} - \E Z^{k+2}.
$$

Of course, if $m=2$, there are no terms on the right-hand side of (3.9) except 
for $\varphi$.

\vskip5mm
\section{{\bf Approximation for $L^r$-norm of Densities $p_n$}}
\setcounter{equation}{0}

\vskip2mm
\noindent
Lemmas 3.2--3.4 can be applied to explore an asymptotic behavior of 
the functionals
$$
I(p) = \|p\|_r^r = \int_{-\infty}^\infty p(x)^r\,dx \qquad (r>1)
$$
with $p = p_n$. Since the densities $p_n$ are well approximated by $\varphi_m$,
we may expect that $I(p_n) \sim I(\varphi_m)$ for large $n$.
However, $\varphi_m$ do not need to be positive on the whole real
line, and it is more natural to consider the integrals
$$
I_T(p) = \int_{|x| \leq T} p(x)^r\,dx, \qquad T>0,
$$
over relatively long intervals. Actually, one may take 
$T = T_n = \sqrt{(s-2)\log n}$ ($s>2$).
By Definition 3.1, for all $n$ large enough,
\be
|\varphi_m(x) - \varphi(x)| \leq \frac{1}{2}\,\varphi(x), \qquad 
|x| \leq T_n,
\en
so $\varphi_m$ is positive on $[-T_n,T_n]$.
On these intervals and for large $n$, consider the functions
$$
\ep_n(x) = \frac{p_n(x) - \varphi_m(x)}{\varphi_m(x)},
$$
so that $p_n(x) = \varphi_m(x) (1 + \ep_n(x))$ with
$|\ep_n(x)| \leq \frac{1}{2}$. Hence, by Taylor's formula, and using (4.1)
together with the non-uniform bound (3.7), we get
\bee
|p_n(x)^r - \varphi_m(x)^r| 
 & \leq & 
c\, \varphi(x)^r\,|\ep_n(x)| \\
& \leq & 
2 c\, \varphi(x)^{r-1}\,|p_n(x) - \varphi_m(x)| \, \leq \,
\delta_n\,\frac{\varphi(x)^{r-1}}{1 + |x|^s}\,n^{-\frac{s-2}{2}}
\ene
with some constant $c$ which does not depend on $x$ and $n \geq n_0$
and some positive sequence $\delta_n \rightarrow 0$. After integration
over $[-T_n,T_n]$, this gives
\be
I_T(p_n) = I_T(\varphi_m) + o(n^{-\frac{s-2}{2}}).
\en
In case $s = m \geq 3$ is integer, by a similar argument based on (3.6), 
we also have
\be
I_T(p_n) = I_T(\varphi_{m-1}) + O(n^{-\frac{s-2}{2}}).
\en

The remaining part of the integral,
$$
J_T(p) = \int_{|x| > T} p(x)^r\,dx,
$$
can be shown to be sufficiently small for $p = p_n$ on the basis of Lemma 3.4.
Indeed, for $Z \sim N(0,1)$, 
$$
\P\{|Z| > T_n\} \leq \frac{1}{T_n}\,e^{-T_n^2/2} = 
o\big(n^{-\frac{s-2}{2}}\big).
$$
Hence, from (4.1) and Definition 3.1, also
$$
\big|\nu_m\{|x| > T_n\}\big| \, = \, o\big(n^{-\frac{s-2}{2}}\big).
$$
Since we assume the smoothness condition (1.2), the densities $p_n$ are 
uniformly bounded by some constant $M$ for all $n \geq n_0$. Therefore,
by Lemma 3.4, for all $n$ large enough,
\bee
J_T(p_n) 
 & \leq &
M^{r-1} \int_{|x| > T_n} p_n(x)\,dx \, = \, M^{r-1} \, \P\{|Z_n| > T_n\} \\
& \leq &
M^{r-1} \, \big|\nu_m\{x: |x| > T_n\}\big| + T_n^{-s}\, o\big(n^{-\frac{s-2}{2}}\big)
\, = \, o\big(n^{-\frac{s-2}{2}}\big).
\ene
Combining this relation with (4.2) and (4.3), we arrive at:

\vskip5mm
{\bf Lemma 4.1.} {\sl Suppose that $\beta_s < \infty$ for $s \geq 2$. 
Then for all $n$ large enough, $Z_n$ have bounded densities $p_n$. Moreover, 
for any $r>1$, as $n \rightarrow \infty$,
\be
\int_{-\infty}^\infty p_n(x)^r\,dx \, = \,
\int_{-T_n}^{T_n} \varphi_m(x)^r\,dx + o\big(n^{-\frac{s-2}{2}}\big), \qquad
m = [s],
\en
where $T_n = \sqrt{(s-2)\log n}$. In particular, if $s = m \geq 3$ is integer, 
we also have
\be
\int_{-\infty}^\infty p_n(x)^r\,dx \, = \,
\int_{-T_n}^{T_n} \varphi_{m-1}(x)^r\,dx + O\big(n^{-\frac{s-2}{2}}\big).
\en
}

\vskip5mm
\section{{\bf Truncated $L^r$-norm of Approximating Densities $\varphi_m$}}
\setcounter{equation}{0}

\vskip2mm
\noindent
Let us now find an explicit expression for the second integral in (4.4),
by applying the Edgeworth approximation
\be
\varphi_m(x) = \varphi(x)\,\bigg(1 + \sum_{k=1}^{m-2} Q_k(x)\,n^{-k/2}\bigg),
\qquad m = [s].
\en

In the case $2 < s < 3$, when $\varphi_m = \varphi_2 = \varphi$,
one may extend the integration in (4.4) to the whole real line
at the expense of the error
$$
\int_{|x| > T} \varphi(x)^r\,dx \, < \, \int_{|x| > T} \varphi(x)\,dx \, = \,
\P\{|Z| > T_n\} \, = \, o\big(n^{-\frac{s-2}{2}}\big),
$$
where 
$
T_n = \sqrt{(s-2)\log n}
$
as before. Hence, (4.4) yields
\be
\int_{-\infty}^\infty p_n(x)^r\,dx \, = \,
\int_{-\infty}^\infty \varphi(x)^r\,dx + o\big(n^{-\frac{s-2}{2}}\big), \qquad
2 < s < 3.
\en
This assertion remains to hold for $s=2$ as well (Theorem 1.1).

Next, assume that $s \geq 3$.
As we know, when $n$ is large enough, $\varphi_m(x)$ is positive for 
$|x| \leq T_n$, so the second integral in (4.4) makes sense, cf. (4.1).
Moreover, in order to raise $\varphi_m(x)$ to the power $r$ on the basis
of (5.1), one may apply the Taylor expansion
\bee
(1 + \ep)^r 
 & = &
1 + r\ep + \frac{r(r-1)}{2!}\,\ep^2 + \dots +
\frac{r(r-1)\dots (r-N+1)}{N!}\,\ep^N + O(\ep^{N+1}) \\
 & = &
1 + \sum_{k = 1}^N \frac{(r)_k}{k!}\,\ep^k + O(\ep^{N+1}), \qquad
N = 1,2,\dots, \ \ep \rightarrow 0,
\ene
where the constant in $O$ depends on $N$ only, as long as 
$|\ep| \leq \frac{1}{2}$. Here we used the standard notation
$(r)_k = r(r-1)\dots (r-k+1)$, with convention $(r)_0 = 1$
to be used later on. Choosing
$$
\ep = \sum_{k=1}^{m-2} Q_k(x)\,n^{-k/2}, \quad |x| \leq T_n,
$$
we have with some constants
depending on the first $m$ absolute moments of $X$ that
$$
|\ep| \, \leq \, \sum_{k=1}^{m-2} |Q_k(x)|\,n^{-k/2} \, \leq \,
C\,(1+|x|)^{3(m-2)} \frac{1}{\sqrt{n}} \, \leq \, 
C'\,\frac{(\log n)^{3(m-2)/2}}{\sqrt{n}} \, \leq \, \frac{1}{2}
$$
for all $n$ large enough in the last inequality. In that case, the above
Taylor expansion is thus valid, i.e., uniformly over all $x \in [-T_n,T_n]$,
as $n \rightarrow \infty$,
\be
(1 + \ep)^r \, = \, 1 + \sum_{k=1}^N \frac{(r)_k}{k!}\,\ep^k + \ep_n(x)
\en
with
$$
\ep_n(x) = O\bigg((1+|x|)^{3(m-2)(N+1)}\frac{1}{n^{(N+1)/2}}\bigg).
$$

Furthermore, by the polynomial formula,
$$
\ep^k \, = \, \sum \frac{k!}{k_1! \dots k_{m-2}!}\
Q_1^{k_1}(x) \dots Q_{m-2}^{k_{m-2}}(x)\,
n^{-\frac{1}{2}\,(k_1 + 2k_2 + \dots + (m-2)\, k_{m-2})},
$$
where the summation is running over all non-negative integers
$k_1,\dots,k_{m-2}$ such that $k_1 + \dots + k_{m-2} = k$. Inserting this
in (5.3) and recalling (5.1), we can represent $\varphi_m(x)^r$ as
$$
\varphi(x)^r
\sum \frac{(r)_{k_1 + \dots + k_{m-2}}}{k_1! \dots k_{m-2}!}\
Q_1^{k_1}(x) \dots Q_{m-2}^{k_{m-2}}(x)\,
n^{-\frac{1}{2}\,(k_1 + 2k_2 + \dots + (m-2)\,k_{m-2})} + \varphi(x)^r \ep_n(x)
$$
with summation over all non-negative integers
$k_1,\dots,k_{m-2}$ such that $k_1 + \dots + k_{m-2} \leq N$.
One may now note that
$$
\int_{-T_n}^{T_n} \varphi(x)^r \ep_n(x) \,dx = 
O\big(n^{-\frac{N+1}{2}}\big).
$$

Let us then choose $N = m-2$. Integrating the above
expression for $\varphi_m(x)^r$ over the interval $[-T_n,T_n]$, 
we can represent $\int_{-T_n}^{T_n} \varphi_m(x)^r\,dx$ as
$$
\sum \frac{(r)_{k_1 + \dots + k_{m-2}}}{k_1! \dots k_{m-2}!}
\int_{-T_n}^{T_n} \varphi(x)^r \,
Q_1^{k_1}(x) \dots Q_{m-2}^{k_{m-2}}(x)\,dx\
\frac{1}{n^{\frac{1}{2}\,(k_1 + 2k_2 + \dots + (m-2)\,k_{m-2})}}
$$
at the expense of an error $O(n^{-\frac{m-1}{2}})$.
Moreover, using the property
$$ 
\int_{|x| \geq T_n} x^N \varphi(x)^r \,dx = o(n^{-\frac{s-2}{2}}),
$$
the above integration may be extended to the whole real line. Hence,
$\int_{-T_n}^{T_n} \varphi_m(x)^r\,dx$ is represented as
$$
\sum \frac{(r)_{k_1 + \dots + k_{m-2}}}{k_1! \dots k_{m-2}!}
\int_{-\infty}^\infty \varphi(x)^r\,
Q_1^{k_1}(x) \dots Q_{m-2}^{k_{m-2}}(x)\,dx\
\frac{1}{n^{\frac{1}{2}\,(k_1 + 2k_2 + \dots + (m-2)\,k_{m-2})}} + 
o\big(n^{-\frac{s-2}{2}}\big).
$$

Here, it is sufficient to keep only the powers of $1/n$
not exceeding $(m-2)/2$. But in that case, for any fixed value of
$$
j = k_1 + 2k_2 + \dots + (m-2)\,k_{m-2},
$$
the constraint $j \leq m-2$ implies that $k_{j+1} = \dots = k_{m-2} = 0$.
That is, for any fixed $j$, we only need to consider the collections
$k_1,\dots,k_j$ of length $j$. Thus, the above representation is
simplified to
\begin{eqnarray}
\int_{-T_n}^{T_n} \varphi_m(x)^r\,dx 
 & = &
\int_{-\infty}^\infty \varphi(x)^r\,dx \nonumber \\
 & & \hskip-30mm + \
\sum \frac{(r)_{k_1 + \dots + k_j}}{k_1! \dots k_j!}
\int_{-\infty}^\infty \varphi(x)^r\,
Q_1^{k_1}(x) \dots Q_j^{k_j}(x)\,dx\,n^{-j/2} + 
o\big(n^{-\frac{s-2}{2}}\big)
\end{eqnarray}
with summation over all $j = 1,\dots,m-2$ and over all non-negative integers
$k_1,\dots,k_j$ such that $k_1 + 2k_2 + \dots + j\,k_j = j$.

As the last simplifying step, we note that $Q_{2k-1}(x)$ represents 
a linear combination of the Hermite polynomials $H_{2i-1}(x)$ and has 
a leading term $x^{3(2k-1)}$ up to a constant. In particular, 
it is an odd function. On the other hand, $Q_{2k}(x)$ represents a linear 
combination of $H_{2i}(x)$'s and has a leading term $x^{6k}$, 
so it is an even function. It follows that any function of the form
\be
Q = Q_1^{k_1}(x) \dots Q_j^{k_j}(x) \qquad (k_1 + 2k_2 + \dots + j\,k_j = j)
\en
is either odd or even, depending on whether $j$ is odd or even. Indeed,
for polynomials of the class 1, defined by
$$
P(x) = c_0 + c_2 x^2 + \dots + c_{2N} x^{2N},
$$
let us put ${\rm Ev}(P) = 2N\, ({\rm mod}\, 2) = 0$, and for the class 2,
defined by 
$$
P(x) = c_1 x + \dots + c_{2N-1}\, x^{2N-1}, 
$$
let us put ${\rm Ev}(P) = 2N-1\, ({\rm mod}\, 2) = 1$. The products of 
such polynomials belong to one of the classes, and we have the property
${\rm Ev}(P_1 P_2) = ({\rm Ev}(P_1) + {\rm Ev}(P_2))\, ({\rm mod}\, 2)$.
Therefore, using ${\rm Ev}(Q_i) = 3i\, ({\rm mod}\, 2) = i\, ({\rm mod}\, 2)$ 
and summation in the group $\mathbb Z_2$, we have
\bee
{\rm Ev}(Q) 
 & = & 
k_1\, {\rm Ev}(Q_1) + \dots + k_j\, {\rm Ev}(Q_j) \\
 & = &
k_1 \cdot 1\, ({\rm mod}\, 2) + \dots + k_j \cdot j\, ({\rm mod}\, 2)
 \, = \,
(k_1 + \dots + j k_j)\, ({\rm mod}\, 2) \, = \, j\, ({\rm mod}\, 2).
\ene

Thus, $Q$ is an odd function in (5.5), as long as $j$ is odd, and then 
the corresponding integral in (5.4) is vanishing. As a result, 
(4.4) and (5.4) yield the following asymptotic expansion, 
which also holds for $2 \leq s < 3$, in view of (5.2).

\vskip5mm
{\bf Proposition 5.1.} {\sl Suppose that $\beta_s < \infty$ for $s \geq 2$. 
Then, with $m =[s]$, for any $r>1$,
\be
\int_{-\infty}^\infty p_n(x)^r\,dx \, = \, 
\int_{-\infty}^\infty \varphi(x)^r\,dx \, \bigg(1 + 
\sum_{j=1}^{[\frac{m-2}{2}]}\, \frac{a_j}{n^j}\bigg) + o\big(n^{-\frac{s-2}{2}}\big)
\en
with coefficients defined by
\be
a_j \int_{-\infty}^\infty \varphi(x)^r\,dx
= \sum \frac{(r)_{k_1 + \dots + k_{2j}}}{k_1! \dots k_{2j}!} 
\int_{-\infty}^\infty 
Q_1^{k_1}(x) \dots Q_{2j}^{k_{2j}}(x)\,\varphi(x)^r\, dx.
\en
Here, the summation is running over all non-negative integers
$k_1,\dots,k_{2j}$ such that $k_1 + 2k_2 + \dots + 2j\,k_{2j} = 2j$, 
with notation $(r)_k = r(r-1)\dots (r-k+1)$.

}

\vskip5mm
From Definition 3.1, it follows that each polynomial $Q_k$
is determined by the moments of $X$ up to order $k+2$. Hence, each $a_j$
in (5.7) is only determined by $r$ and by the moments -- or equivalently,
by the cumulants of $X$ up to order $2j+2$. Moreover, $a_j = 0$ if
these cumulants are vanishing.

\vskip5mm
\section{{\bf The Case where the First Cumulants are Vanishing}}
\setcounter{equation}{0}

\vskip2mm
\noindent
For $2 \leq s < 4$, we necessarily have $m \leq 3$, so that 
the sum in (5.6) has no term, and then
\be
\int_{-\infty}^\infty p_n(x)^r\,dx \, = \,
\int_{-\infty}^\infty \varphi(x)^r\,dx + o\big(n^{-\frac{s-2}{2}}\big).
\en

In the more interesting case $s \geq 4$, the leading term in the Edgeworth 
expansion (5.1) may be written explicitly, as was already done in
the representation (3.9). It implies that, for some unique $1 \leq k \leq m-2$,
\be
\varphi_m(x) \, = \,
\varphi(x) + \varphi(x) \frac{\gamma_{k+2}}{(k+2)!}\, H_{k+2}(x)\,n^{-k/2} +
C(x) \varphi(x)\,(1 + |x|^{3(m-2)})\,n^{-(k+1)/2}
\en
with some function $C(x)$ bounded by a constant which does not depend
on $x$ and large $n \geq n_0$. To study an asymptotic behavior of the
truncated $L^r$-norm of $\varphi_m$, one may repeat computations of the previous
section in this simple particular case, or alternatively, one may
just refer to the general result described in Proposition 5.1.

Indeed, (6.2) is equivalent to saying that the first moments of $X$ up to order 
$k+1$ coincide with those of $Z \sim N(0,1)$ for some $1 \leq k \leq m-2$.
Therefore, as emphasized after Proposition 5.1,
$a_j = 0$ whenever $2j+2 \leq k+1$, that is, $j \leq \frac{k-1}{2}$.
Then also $Q_j = 0$. In case $2j+2 = k+2$, that is,
$j = k/2$ with even $k$, all terms in the sum (5.7) are vanishing, except
(potentially) for the term corresponding to
$k_1 = \dots = k_{2j-1} = 0$, $k_{2j} = 1$. Then the right-hand side of
(5.7) becomes
$$
r \int_{-\infty}^\infty Q_{2j}(x)\,\varphi(x)^r\, dx =
r \int_{-\infty}^\infty Q_k(x)\,\varphi(x)^r\, dx = 
r\, \frac{\gamma_{k+2}}{(k+2)!}
\int_{-\infty}^\infty H_{k+2}(x)\,\varphi(x)^r\, dx,
$$
and hence (5.6) yields
\be
\int_{-\infty}^\infty p_n(x)^r\,dx \, = \,
\int_{-\infty}^\infty \varphi(x)^r\,dx + A n^{-k/2} +
O(n^{-\frac{k+1}{2}}) + o\big(n^{-\frac{s-2}{2}}\big),
\en
where
$$
A = r\,\frac{\gamma_{k+2}}{(k+2)!}\,
\int_{-\infty}^\infty H_{k+2}(x)\,\varphi(x)^r\, dx, \qquad
\gamma_{k+2} = \E X^{k+2} - \E Z^{k+2}.
$$
In particular, $A = 0$ whenever $k$ is odd (since the corresponding
Chebyshev-Hermite polynomial is odd).

To proceed, we need to focus on the integrals of the form
$I(k,r) = \int_{-\infty}^\infty H_k(x)\,\varphi(x)^r\,dx$ with even $k$.

\vskip5mm
{\bf Lemma 6.1.} {\sl For any $k = 1,2,\dots$,
\be
I(2k,r) \, = \,
\frac{(2k-1)!!}{r^{\frac{2k+1}{2}}\, (2\pi)^{\frac{r-1}{2}}}\ (1-r)^k.
\en
}

{\bf Proof.} The $k$-th Chebyshev-Hermite polynomial
\be
H_k(x) \, = \, (-1)^k\, \big(e^{-x^2/2}\big)^{(k)}\, e^{x^2/2} \, = \,
\E\,(x + iZ)^k, \qquad Z \sim N(0,1),
\en
has generating function
$$
\sum_{k=0}^\infty H_k(x)\, \frac{z^k}{k!} = e^{xz - z^2/2}, \qquad z \in \C,
$$
from which one can find the generating function for the sequence
$c_k = I(k,r)$. Namely,
$$
\sum_{k=0}^\infty c_k\, \frac{z^k}{k!} \, = \,
\int_{-\infty}^\infty e^{xz - z^2/2}\,\varphi(x)^r\,dx \, = \,
\frac{1}{(2\pi)^{\frac{r-1}{2}}\sqrt{r}}\,e^{-\frac{1}{2}\,(1 - \frac{1}{r})\,z^2}.
$$
Differentiating this equality $2k$ times and applying the definition (6.5),
we arrive at
$$
c_{2k} = \frac{1}{(2\pi)^{\frac{r-1}{2}}\sqrt{r}}\,
\Big(1 - \frac{1}{r}\Big)^k\,H_{2k}(0).
$$
It remains to apply the second equality in (6.5), which gives
$$ H_{2k}(0) = (-1)^k\,\E Z^{2k} = (-1)^k\,(2k-1)!! $$
\qed

For the first three even values $k = 2,4,6$, we thus have
\begin{eqnarray}
I(2,r) & = & -\frac{1}{r^{3/2}\,(2\pi)^{\frac{r-1}{2}}}\, (r-1), \qquad
I(4,r) \ = \ \frac{3}{r^{5/2}\,(2\pi)^{\frac{r-1}{2}}}\,(r-1)^2, \nonumber \\
I(6,r) & = & -\frac{15}{r^{7/2}\,(2\pi)^{\frac{r-1}{2}}}\,(r-1)^3.
\end{eqnarray}

With similar arguments, one may also evaluate the integrals
$\int_{-\infty}^\infty H_k(x)^2\,\varphi(x)^r\,dx$. For example,
\be
\int_{-\infty}^\infty H_3(x)^2\,\varphi(x)^r\,dx
 \, = \,
\frac{1}{\sqrt{r}\,(2\pi)^{\frac{r-1}{2}}}\
\E\,
\bigg(\Big(\frac{Z}{\sqrt{r}}\Big)^3 - 3\, \Big(\frac{Z}{\sqrt{r}}\Big)\bigg)^2 
 \, = \,
\frac{3\, (5 - 6r + 3r^2)}{r^{7/2}\,(2\pi)^{\frac{r-1}{2}}}.
\en

Thus, the formula (6.4) may be used in the asymptotic representation (6.3).
The particular case $k = [s]-2$ should be mentioned separately.

\vskip5mm
{\bf Corollary 6.2.} {\sl Suppose that $\E X^l = \E Z^l$ for
$l = 1,\dots,m-1$ $(m \geq 3)$, where $Z \sim N(0,1)$.
If $\beta_s < \infty$ for some $s \in [m,m+1)$, then for all $n$ large enough, 
$Z_n$ have bounded densities $p_n$. Moreover,
\be
\int_{-\infty}^\infty p_n(x)^r\,dx \, = \,
\int_{-\infty}^\infty \varphi(x)^r\,dx + A n^{-\frac{m-2}{2}} +
o\big(n^{-\frac{s-2}{2}}\big)
\en
with $A = 0$ in the case $m = 2k-1$ is odd, while in the case where 
$m = 2k$ is even, we have
$$
A = \frac{\gamma_{2k}}{2^k k!}\,
\frac{1}{(2\pi)^{\frac{r-1}{2}}}\, \frac{(1-r)^k}{r^{\frac{2k-1}{2}}}, \qquad
\gamma_{2k} = \E X^{2k} - \E Z^{2k}.
$$
If $\beta_s < \infty$ for $s = m+1$, then $o$-term 
in $(6.8)$ may be replaced with $O$-term.
}

\vskip5mm
For example, if $\gamma_3 = \E X^3 = 0$, so that $m=4$, $4 \leq s < 5$, we have
$$
A = \frac{\gamma_4}{8}\,
\frac{1}{(2\pi)^{\frac{r-1}{2}}}\, \frac{(1-r)^2}{r^{\frac{3}{2}}}, \qquad
\gamma_4 = \E X^4 - \E Z^4 = \E X^4 - 3,
$$
and (6.8) becomes
\be
\int_{-\infty}^\infty p_n(x)^r\,dx \, = \,
\int_{-\infty}^\infty \varphi(x)^r\,dx + A n^{-1} +
o\big(n^{-\frac{s-2}{2}}\big).
\en
By (6.3), a similar formula remains to hold in the case
$5 \leq s < 6$, but then the $o$-term should be replaced with
$O(n^{-3/2})$.

\vskip5mm
\section{{\bf Moments of Order $4 \leq s \leq 8$}}
\setcounter{equation}{0}

\vskip2mm
\noindent
Returning to the general expansion (5.6) in Proposition 5.1 with
coefficients $a_j$ described in (5.7), let us now derive formulas similar 
to (6.9) for two regions of the values of $s$ without additional assumptions 
on the first cumulants. To evaluate the integrals in that definition, we
will use the formulas for the polynomials $Q_j$ described in Section 3
for the indexes $j \leq 4$.

If $4 \leq s < 6$, the expansion (5.6) contains only one term, namely, we get 
\be
\int_{-\infty}^\infty p_n(x)^r\,dx \, = \, 
\int_{-\infty}^\infty \varphi(x)^r\,dx + \frac{a_1}{n}
\int_{-\infty}^\infty \varphi(x)^r\,dx  + o\big(n^{-\frac{s-2}{2}}\big)
\en
with the coefficient for $j=1$ in front of $1/n$, i.e.,
\bee
A_1 \, \equiv \, a_1 \int_{-\infty}^\infty \varphi(x)^r\,dx
 & = &
\frac{(r)_1}{1!} \int_{-\infty}^\infty Q_2(x)\,\varphi(x)^r\, dx +
\frac{(r)_2}{2!} \int_{-\infty}^\infty Q_1^2(x)\,\varphi(x)^r\, dx
 \\
 & = & 
r \int_{-\infty}^\infty \Big(\frac{\gamma_4}{4!}\, H_4(x) + 
\frac{1}{2!} \left(\frac{\gamma_3}{3!}\right)^2 H_6(x)\Big)\,\varphi(x)^r\, dx \\
 & & + \
\frac{r(r-1)}{2} \int_{-\infty}^\infty 
\Big(\frac{\gamma_3}{3!}\, H_3(x)\Big)^2\,\varphi(x)^r\, dx.
\ene
Applying the formulas (6.6)-(6.7), we find that
\bee
A_1 
 & = &
r\,\frac{\gamma_3^2}{2!\,3!^2}\, I(6,r) + 
r\,\frac{\gamma_4}{4!}\,I(4,r) +
\frac{r(r-1)}{2}\,\Big(\frac{\gamma_3}{3!}\Big)^2 
\int_{-\infty}^\infty H_3(x)^2\, \varphi(x)^r\,dx \\
 & = &
r\,\frac{\gamma_3^2}{72}\, I(6,r) + 
r\,\frac{\gamma_4}{24}\,I(4,r) +
r(r-1)\,\frac{\gamma_3^2}{72} \,
\frac{3}{r^{7/2}\,(2\pi)^{\frac{r-1}{2}}}\,(5 - 6r + 3r^2)
\\ 
 & = &
-r\,\frac{\gamma_3^2}{72}\,\frac{15}{r^{7/2}\,(2\pi)^{\frac{r-1}{2}}}\,(r-1)^3 + 
r\,\frac{\gamma_4}{24}\,\frac{3}{r^{5/2}\,(2\pi)^{\frac{r-1}{2}}}\,(r-1)^2 \\
 & & + \ 
r(r-1)\,\frac{\gamma_3^2}{24}\,
\frac{1}{r^{7/2}\,(2\pi)^{\frac{r-1}{2}}}\,(5 - 6r + 3r^2).
\ene
Equivalently,
$$
(2\pi)^{\frac{r-1}{2}} \frac{r^{5/2}}{r-1}\ A_1 = 
-\frac{5}{24}\,(r-1)^2\, \gamma_3^2 + 
\frac{1}{8}\,r(r-1)\, \gamma_4 + \frac{1}{24}\,
(5 - 6r + 3r^2)\,\gamma_3^2.
$$
Collecting the coefficients in front of $\gamma_3^2$, we arrive at the
following refinement of (7.1).

\vskip5mm
{\bf Proposition 7.1.} {\sl Suppose that $\beta_s < \infty$ for $4 \leq s < 6$. Then, for any $r>1$, 
\be
\int_{-\infty}^\infty p_n(x)^r\,dx \, = \,
\int_{-\infty}^\infty \varphi(x)^r\,dx + A_1 n^{-1} + 
o\big(n^{-\frac{s-2}{2}}\big),
\en
where the constant $A_1 = A_1(r)$ is given by
\be
(2\pi)^{\frac{r-1}{2}} \frac{r^{3/2}}{r-1}\ A_1(r) \, =  \,
\frac{2-r}{12}\, \gamma_3^2 + \frac{r-1}{8}\, \gamma_4.
\en
In the case $s=6$, the formula $(7.2)$ remains valid with the remainder term
$O(n^{-2})$.
}

\vskip5mm
Note that
\be
\lim_{r \rightarrow 1}\ \frac{A_1(r)}{r-1} \, = \, \frac{1}{12}\,\gamma_3^2.
\en

If $\gamma_3 = 0$, then (7.3) is simplified to
$$
(2\pi)^{\frac{r-1}{2}} \frac{r^{3/2}}{r-1}\ A_1(r) \, =  \, 
\frac{r-1}{8}\, \gamma_4,
$$
which is exactly the constant $A$ in the equality (6.9),
obtained under the cumulant conditions.

Let us now consider the region $6 \leq s < 8$. In this case, the sum in (5.6) 
contains two terms, proportional to $\frac{1}{n}$ and $\frac{1}{n^2}$. 
The coefficient $a_1$ will be as before, while according to (5.7),
\bee
a_2  \int_{-\infty}^\infty \varphi(x)^r\,dx
 & = &
r \int_{-\infty}^\infty Q_4(x)\,\varphi(x)^r\, dx \\
 & & + \ 
\frac{(r)_2}{2} \int_{-\infty}^\infty 
(Q_2^2(x) + 2Q_1(x)Q_3(x))\,\varphi(x)^r\, dx \\
 & & + \ 
\frac{(r)_3}{2} \int_{-\infty}^\infty Q_1^2(x) Q_2(x)\,\varphi(x)^r\, dx +
\frac{(r)_4}{24} \int_{-\infty}^\infty Q_1^4(x)\,\varphi(x)^r\, dx.
\ene
We thus have the following refinement of Proposition 7.1 under stronger
moment assumptions.

\vskip5mm
{\bf Proposition 7.2.} {\sl Suppose that $\beta_s < \infty$ for $6 \leq s < 8$. Then, for any $r>1$, 
\be
\int_{-\infty}^\infty p_n(x)^r\,dx \, = \,
\int_{-\infty}^\infty \varphi(x)^r\,dx + A_1 n^{-1} + A_2 n^{-2} + 
o\big(n^{-\frac{s-2}{2}}\big),
\en
where $A_1$ is given in $(7.3)$ and
\bee
A_2 
 & = &
r \int_{-\infty}^\infty Q_4(x)\,\varphi(x)^r\, dx +
\frac{(r)_2}{2} \int_{-\infty}^\infty 
\big(Q_2^2(x) + 2\,Q_1(x)Q_3(x)\big)\,\varphi(x)^r\, dx \\
 & & + \
\frac{(r)_3}{2} \int_{-\infty}^\infty Q_1^2(x) Q_2(x)\,\varphi(x)^r\, dx +
\frac{(r)_4}{24} \int_{-\infty}^\infty Q_1^4(x)\,\varphi(x)^r\, dx.
\ene
In the case $s=6$, the formula $(7.4)$ remains valid with the remainder term
$O(n^{-3})$.
}

\vskip5mm
We can rewrite $A_2$ in terms of the cumulants of $X$ as follows:
\begin{eqnarray*}
A_2 
 & = &
r \frac{\gamma_6}{6!} 
\int_{-\infty}^\infty H_6(x)\, \varphi(x)^r\,dx + 
r \Big(\frac{\gamma_3 \gamma_5}{3!\,5!} + \frac{\gamma_4^2}{2!\,4!^2}\Big) 
\int_{-\infty}^\infty H_8(x)\, \varphi(x)^r\,dx \\ 
 & & + \
r \, 
\frac{\gamma_3^2 \gamma_4}{2!\,3!^2\,4!} 
\int_{-\infty}^\infty H_{10}(x)\, \varphi(x)^r\,dx + 
r \, 
\frac{\gamma_3^4}{4!\,3!^4} \int_{-\infty}^\infty H_{12}(x)\, \varphi(x)^r\,dx \\
 & & + \ 
\frac{r(r-1)}{2}\, 
\frac{\gamma_4^2}{4!^2} \int_{-\infty}^\infty H_4(x)^2 \varphi(x)^r\,dx + 
\frac{r(r-1)}{2}\,\frac{\gamma_3^4}{2!^2\,3!^4} 
\int_{-\infty}^\infty H_6(x)^2 \varphi(x)^r\,dx \\ 
 & & + \ 
\frac{r(r-1)}{2}\, 
\frac{\gamma_3^2 \gamma_4}{3!^2\,4!}
\int_{-\infty}^\infty H_4(x) H_6(x)\, \varphi(x)^r\,dx +
r(r-1)\, 
\frac{\gamma_3\gamma_5}{3!\,5!} 
\int_{-\infty}^\infty H_3(x) H_5(x)\, \varphi(x)^r\,dx \\
 & & + \
r(r-1)\,
\frac{\gamma_3^2 \gamma_4}{3!^2\,4!} 
\int_{-\infty}^\infty H_3(x) H_7(x)\, \varphi(x)^r\,dx +
r(r-1)\,
\frac{\gamma_3^4}{(3!)^5} 
\int_{-\infty}^\infty H_3(x) H_9(x)\, \varphi(x)^r\,dx \\ 
 & & + \ 
\frac{r(r-1)(r-2)}{2}\, 
\frac{\gamma_3^2 \gamma_4}{3!^2\,4!} 
\int_{-\infty}^\infty H_3(x)^2 H_4(x)\, \varphi(x)^r\,dx \\ 
 & & + \
\frac{r(r-1)(r-2)}{4}\, 
\frac{\gamma_3^4}{3!^4} 
\int_{-\infty}^\infty H_3(x)^2 H_6(x)\, \varphi(x)^r\,dx \\ 
 & & + \ 
\frac{r(r-1)(r-2)(r-3)}{24}\, \frac{\gamma_3^4}{3!^4} 
\int_{-\infty}^\infty H_3(x)^4 \varphi(x)^r\,dx.
\end{eqnarray*}

In the case $\gamma_3 = 0$, this long expression is simplified to
\bee
A_2 
 & = &
r \frac{\gamma_6}{6!} \int_{-\infty}^\infty H_6(x)\, \varphi(x)^r\,dx \\
 & & + 
r \frac{\gamma_4^2}{2!\,4!^2} \int_{-\infty}^\infty H_8(x)\, \varphi(x)^r\,dx +
r(r-1)\, 
\frac{\gamma_4^2}{2!\,4!^2} \int_{-\infty}^\infty H_4(x)^2\, \varphi(x)^r\,dx. 
\ene

\vskip5mm
\section{{\bf Expansions for R\'enyi Entropies}}
\setcounter{equation}{0}

\vskip2mm
\noindent
Let us now reformulate the asymptotic results about 
the integrals $\int_{-\infty}^\infty p_n(x)^r\,dx$
in terms of R\'enyi's entropies and entropy powers
$$
h_r(Z_n) \, = \, -\frac{1}{r-1}\,\log \int_{-\infty}^\infty p_n(x)^r\,dx, \qquad
N_r(Z_n) \, = \, 
\bigg(\int_{-\infty}^\infty p_n(x)^r\,dx\bigg)^{-\frac{2}{r-1}}.
$$
Since these functionals represent smooth functions of the $L^r$-norm,
from Proposition 5.1 we immediately obtain:

\vskip5mm
{\bf Proposition 8.1.} {\sl Let $\E\,|X|^s < \infty$ for some 
$s \geq 2$, and $m =[s]$. Then, for any $r>1$,
\begin{eqnarray}
h_r(Z_n) 
 & = & h_r(Z) + \sum_{j=1}^{[\frac{m-2}{2}]}\, 
\frac{b_j}{n^j} + o\big(n^{-\frac{s-2}{2}}\big), \\
N_r(Z_n) 
 & = & N_r(Z)\,\bigg(1 +  
\sum_{j=1}^{[\frac{m-2}{2}]}\, \frac{c_j}{n^j}\bigg) + 
o\big(n^{-\frac{s-2}{2}}\big)
\end{eqnarray}
with coefficients $b_j$ and $c_j$ that are determined by $r$ and 
by the moments of $X$ up to order~$2j+2$.
}

\vskip5mm
{\bf Proof of Theorem 1.2.}
To evaluate the first coefficients in the expansions (8.1)-(8.2), 
we apply Taylor's formulas
\begin{eqnarray}
\log(a + b + c) 
 & = &
\log a + a^{-1} b + O(b^2 + |c|), \\
(a + b + c)^q 
 & = & 
a^q + q a^{q-1}\, b + O(b^2 + |c|), \nonumber
\end{eqnarray}
holding with $a>0$, $q \neq 0$, and $b,c \rightarrow 0$.
For $q = -\frac{2}{r-1}$, the last equality reads
\be
(a + b + c)^{-\frac{2}{r-1}} \, = \, 
a^{-\frac{2}{r-1}} -\frac{2}{r-1}\, a^{-\frac{r+1}{r-1}}\, b + O(b^2 + |c|).
\en
In particular (with $b=0$), the expansion of the form
$$
\int_{-\infty}^\infty p_n(x)^r\,dx \, = \,
\int_{-\infty}^\infty \varphi(x)^r\,dx + o\big(n^{-\frac{s-2}{2}}\big),
$$
which corresponds in Proposition 5.1 to the region $2 < s < 4$, implies
$$
\log \int_{-\infty}^\infty p_n(x)^r\,dx =
\log \int_{-\infty}^\infty \varphi(x)^r\,dx + o\big(n^{-\frac{s-2}{2}}\big).
$$
Equivalently, $h_r(Z_n) = h_r(Z) + o(n^{-\frac{s-2}{2}})$ or
$N_r(Z_n) = N_r(Z) + o(n^{-\frac{s-2}{2}})$ for $Z \sim N(0,1)$. 

More generally, applying (8.3)-(8.4) to the expansion
$$
\int_{-\infty}^\infty p_n(x)^r\,dx \, = \,
\int_{-\infty}^\infty \varphi(x)^r\,dx + 
A_1\, n^{-1} + o\big(n^{-\frac{s-2}{2}}\big),
$$
corresponding to Proposition 7.1 with its region $4 \leq s < 6$, we get
$$
\log \int_{-\infty}^\infty p_n(x)^r\,dx \ = \
\log \int_{-\infty}^\infty \varphi(x)^r\,dx +
A_1\, n^{-1}\, \bigg(\int_{-\infty}^\infty \varphi(x)^r\,dx\bigg)^{-1} +
o\big(n^{-\frac{s-2}{2}}\big),
$$
and
\bee
\bigg(\int_{-\infty}^\infty p_n(x)^r\,dx\bigg)^{-\frac{2}{r-1}}
 & = &
\bigg(\int_{-\infty}^\infty \varphi(x)^r\,dx\bigg)^{-\frac{2}{r-1}} \\
 & &
-\frac{2A_1}{r-1}\, n^{-1}\,
\bigg(\int_{-\infty}^\infty \varphi(x)^r\,dx\bigg)^{-\frac{r+1}{r-1}} +
o\big(n^{-\frac{s-2}{2}}\big).
\ene
Thus,
\be
h_r(Z_n) = h_r(Z) - \frac{A_1}{r-1}\,N_r(Z)^{\frac{r-1}{2}}\, n^{-1}
+ o(n^{-\frac{s-2}{2}}),
\en
and (equivalently)
\begin{eqnarray}
N_r(Z_n) 
 & = &
N_r(Z) -\frac{2A_1}{r-1}\,N_r(Z)^{\frac{r+1}{2}}\, n^{-1} + 
o(n^{-\frac{s-2}{2}}) \nonumber \\
 & = &
N_r(Z)\, \Big[1 -\frac{2A_1}{r-1}\,N_r(Z)^{\frac{r-1}{2}}\, n^{-1}\Big] + 
o(n^{-\frac{s-2}{2}}).
\end{eqnarray}

Recall that $A_1 = A_1(r)$ is determined by $r$ and the cumulants 
$\gamma_3 = \E X^3$ and $\gamma_4 = \E X^4 - 3$. More precisely, according
to the formula (7.3) of Proposition 7.1,
$$
\frac{A_1}{r-1} \, =  \,
\frac{1}{(2\pi)^{\frac{r-1}{2}}\,r^{3/2}}\ \bigg[
\frac{2-r}{12}\, \gamma_3^2 + \frac{r-1}{8}\, \gamma_4\bigg].
$$
Since also
$$
N_r(Z)^{\frac{r-1}{2}} = 
\bigg(\int_{-\infty}^\infty \varphi(x)^r\,dx\bigg)^{-1} =
(2\pi)^{\frac{r-1}{2}}\, r^{1/2},
$$
the coefficients $b_1$ and $c_1$ in (8.1)-(8.2) in front of $n^{-1}$
are simplified according to (8.5)-(8.6) as
$$
b_1 = -\frac{A_1}{r-1}\,N_r(Z)^{\frac{r-1}{2}} = -\frac{1}{r}\ \bigg[
\frac{2-r}{12}\, \gamma_3^2 + \frac{r-1}{8}\, \gamma_4\bigg], \qquad c_1 = 2b_1.
$$
\qed

\vskip2mm
Let us complement the expansions of Theorem 1.2 with similar assertions
corresponding to the scenario from Corollary 6.2, where
the first $m-1$ moments of $X$ coincide with those of $Z \sim N(0,1)$,
for some integer $m \geq 3$. If $\beta_s$ is finite for $s \in [m,m+1)$,
in that case we have an expansion of the form
$$
\int_{-\infty}^\infty p_n(x)^r\,dx \, = \,
\int_{-\infty}^\infty \varphi(x)^r\,dx + 
A n^{-\frac{m-2}{2}} + o\big(n^{-\frac{s-2}{2}}\big).
$$
Hence, by (8.3)-(8.4),
\bee
\log \int_{-\infty}^\infty p_n(x)^r\,dx
 & = &
\log \int_{-\infty}^\infty \varphi(x)^r\,dx \\
 & & + \ 
A\, n^{-\frac{m-2}{2}}\,
\bigg(\int_{-\infty}^\infty \varphi(x)^r\,dx\bigg)^{-1} +
O(n^{-(m-2)}) + o\big(n^{-\frac{s-2}{2}}\big),
\ene
and
\bee
\bigg(\int_{-\infty}^\infty p_n(x)^r\,dx\bigg)^{-\frac{2}{r-1}}
 & = &
\bigg(\int_{-\infty}^\infty \varphi(x)^r\,dx\bigg)^{-\frac{2}{r-1}} \\
 & & \hskip-10mm
-\frac{2A}{r-1}\, n^{-\frac{m-2}{2}}\,
\bigg(\int_{-\infty}^\infty \varphi(x)^r\,dx\bigg)^{-\frac{r+1}{r-1}} +
O(n^{-(m-2)}) + o\big(n^{-\frac{s-2}{2}}\big).
\ene
Since $m-2 > \frac{s-2}{2}$, here $O$-term may be removed. 
In addition, as before, the last integral with its
power can be written as $N_r(Z)^{\frac{r+1}{2}}$.
Therefore, we obtain the asymptotic relations
$$
h_r(Z_n) = 
h_r(Z) - \frac{A}{r-1}\,N_r(Z)^{\frac{r-1}{2}}\, n^{-\frac{m-2}{2}}
+ o(n^{-\frac{s-2}{2}})
$$
and
\bee
N_r(Z_n) 
 & = &
N_r(Z) - \frac{2A}{r-1}\,N_r(Z)^{\frac{r+1}{2}}\, n^{-\frac{m-2}{2}}
+ o(n^{-\frac{s-2}{2}}) \\ 
 & = &
N_r(Z)\,\bigg[1 -
\frac{2A}{r-1}\,N_r(Z)^{\frac{r-1}{2}}\, n^{-\frac{m-2}{2}}\bigg]
+ o(n^{-\frac{s-2}{2}})
\ene
in full analogy with (8.5)-(8.6). The only difference is that
we have a different formula for the constant $A = A(r)$.
As stated in Corollary 6.2, here $A = 0$ in the case $m = 2k-1$ is odd, 
while in the case $m = 2k$ is even, we have
$$
A = \frac{\gamma_{2k}}{2^k k!}\,
\frac{1}{(2\pi)^{\frac{r-1}{2}}}\, \frac{(1-r)^k}{r^{\frac{2k-1}{2}}}, \qquad
\gamma_{2k} = \E X^{2k} - \E Z^{2k}.
$$
Using again
$
N_r(Z)^{\frac{r-1}{2}} = 
(2\pi)^{\frac{r-1}{2}}\, r^{1/2},
$
the coefficients $b_{k-1}$ and $c_{k-1}$ in (8.1)-(8.2) in front of 
$n^{-\frac{m-2}{2}} = n^{-(k-1)}$ are simplified to
$$
b_{k-1} = -\frac{A}{r-1}\,N_r(Z)^{\frac{r-1}{2}} = 
\frac{\gamma_{2k}}{2^k k!}\,\frac{(1-r)^{k-1}}{r^{k-1}}, \qquad
c_{k-1} = 2b_{k-1}.
$$
Let us also remind that, if $\beta_s < \infty$ for $s = m+1$, then $o$-term 
may be replaced with $O(n^{-\frac{m-1}{2}})$.
We are thus ready to make a corresponding statement.

\vskip5mm
{\bf Proposition 8.2.} {\sl Suppose that
$\E X^l = \E Z^l$ for $l = 3,\dots,m-1$ $(m \geq 3)$.
If $\beta_s < \infty$ for some $s \in [m,m+1)$, then for any $r>1$,
\bee
h_r(Z_n) 
 & = &
h_r(Z) + b n^{-\frac{m-2}{2}} + o(n^{-\frac{s-2}{2}}), \\
N_r(Z_n) 
 & = &
N_r(Z)\,\big(1 + 2 b\, n^{-\frac{m-2}{2}}\big) + o(n^{-\frac{s-2}{2}})
\ene
with constant $b = 0$ in the case $m = 2k-1$ is odd, 
while in the case $m = 2k$ is even,
$$
b = b_{k-1} = \frac{\gamma_{2k}}{2^k k!}\,\Big(\frac{1}{r} - 1\Big)^{k-1}, 
\qquad \gamma_{2k} = \E X^{2k} - \E Z^{2k}.
$$
If $\beta_s < \infty$ for $s = m+1$, then $o$-term may be replaced
with $O(n^{-\frac{m-1}{2}})$.
}

\vskip5mm
For example, if $\gamma_3 = \E X^3 = 0$, we return to the equality (1.4) 
from Theorem 1.2.

\vskip5mm
\section{{\bf Comparison with the entropic CLT. Monotonicity}}
\setcounter{equation}{0}

\vskip2mm
\noindent
Put
$$
\Delta_n(r) = h_r(Z) - h_r(Z_n), \qquad \Delta_n = \Delta_n(1).
$$
The latter quantity, which may also be written as 
$D(Z_n||Z) = \int_{-\infty}^\infty p_n(x)\,\log\frac{p_n(x)}{\varphi(x)}\,dx$,
represents the Kullback-Leibler distance from the distribution of $Z_n$
to the standard normal law (or, the relative entropy).
As was mentioned, the sequence $\Delta_n$ is always non-negative
and non-increasing. Moreover, the entropic CLT
asserts that $\Delta_n \rightarrow 0$ as $n \rightarrow \infty$,
as long as $\Delta_n$ is finite for some $n$
(in general, it is a weaker condition in comparison with (1.2)).
The basic references for these results are [Ba], [A-B-B-N], [M-B].

The rate of convergence of $\Delta_n$ to zero was studied in [B-C-G2],
and here we recall a few asymptotic results, assuming that $\Delta_n < \infty$
for some $n$, and that $\beta_s = \E\,|X|^s < \infty$ for a real number
$s \geq 2$. Namely, we have
$$
\Delta_n = o\bigg(\frac{1}{(n \log n)^{\frac{s-2}{2}}}\bigg), \qquad 2 \leq s < 4.
$$
Modulo a logarithmic term, it is the same rate as for $\Delta_n(r)$
indicated in Theorem 1.2. Nevertheless, it is not yet clear, if
one can similarly improve Theorem 1.2. On the other hand, 
for any prescribed $\eta>1$, it may occur that, for all $n$ large enough, 
$$
\Delta_n \geq \frac{c}{(n \log n)^{\frac{s-2}{2}}\,(\log n)^\eta}
$$
with some constant $c = c(\eta,s)>0$ depending on $\eta$ and $s$ only 
([B-C-G2], Theorem 1.3).

The range $s \geq 4$ is more interesting, since then one may control
the speed of $\Delta_n$. In particular,
\bee
\Delta_n 
 & = & 
\frac{\gamma_3^2}{12}\,n^{-1} +
o\bigg(\frac{1}{(n \log n)^{\frac{s-2}{2}}}\bigg), \qquad 4 \leq s < 6, \\
\Delta_n 
 & = &
\frac{\gamma_3^2}{12}\,n^{-1} +
O\bigg(\frac{1}{(n \log n)^2}\bigg), \qquad \quad s = 6.
\ene
Thus, if $\gamma_3 \neq 0$, then $\Delta_n$ is equivalent to a decreasing
sequence, which decreases at rate $n^{-1}$. (Strictly speaking, 
this property does not imply the monotonicity itself.)

Let us compare this asymptotic with what is given in Theorem 1.2.
Namely, for any $r>1$, we have
\begin{eqnarray}
\Delta_n(r)
 & = & 
B_1\, n^{-1} + o\big(n^{-\frac{s-2}{2}}\big), \qquad 4 \leq s < 6, \\
\Delta_n(r)
 & = &
B_1\, n^{-1} + O\big(n^{-2}\big), \qquad \quad s = 6,
\end{eqnarray}
where
$$
B_1 = B_1(r) = - b = \frac{1}{4r}\, \bigg[
\frac{2-r}{3}\, \gamma_3^2 + \frac{r-1}{2}\, \gamma_4\bigg].
$$
We see that $B(r) \rightarrow \frac{1}{12}\,\gamma_3^2$ as $r \rightarrow 1$,
so that we recover the main term in the asymptotic for $\Delta_n$,
and at the same rate modulo a logarithmic factor. 

However, what can one say about the sign of $B_1(r)$ with fixed $r>1$?
First suppose that $\gamma_3 \neq 0$. When $r$ is sufficiently close to 1,
then $B_1(r)>0$, so that $\Delta_n(r)$ is equivalent to a decreasing 
sequence like for $r=1$. More precisely, this is true for all $r>1$,
whenever $\gamma_4 \geq \frac{2}{3}\,\gamma_3^2$. But, if 
$\gamma_4 < \frac{2}{3}\,\gamma_3^2$, then $B_1(r) < 0$ for all 
$$
r > r_0 = \frac{4\gamma_3^2 - 3\gamma_4}{2\gamma_3^2 - 3\gamma_4}.
$$ 
Hence  $\Delta_n(r)$ becomes to be equivalent to 
an increasing sequence. In that case, necessarily
$h_r(Z_n) > h_r(Z)$ for all $n$ large enough,
which is impossible in the Shannon case $r=1$. This shows that
$\Delta_n(r)$ may not serve as~distance!

If $\gamma_3 = 0$ (as in case of symmetric distributions), the constant
is simplified to
$$
B_1 = B_1(r) = \frac{r-1}{8r}\, \gamma_4, \qquad \gamma_4 = \E X^4 - 3,
$$
and then the sign of $B_1$ coincides with the sign of $\gamma_4$.
Both cases, $\gamma_4>0$ or $\gamma_4<0$, are typical, and one
can make a similar conclusion as before, but for the whole range $r>1$.
Namely, if $\gamma_4>0$, then $\Delta_n(r)$ is equivalent 
to a decreasing sequence, which decreases at rate $n^{-1}$, and if 
$\gamma_4<0$, then $\Delta_n(r)$ is equivalent 
to an increasing sequence, which increases also at rate $n^{-1}$.

In order to make a more rigorous conclusion about the monotonicity
of $\Delta_n(r)$ for large $n$, the expansions for Renyi entropy
$h_r(Z_n)$ such as (9.1)-(9.2), 
are insufficient. We need to use more terms in the general Proposition 8.1 
involving the quadratic terms $b_2/n^2$ and $c_2/n^2$. This is possible 
under stronger moment assumptions, corresponding to the range $6 \leq s < 8$. 
Indeed, in that case, Proposition 8.1 provides the expansion (1.5)
in which the coefficient $b_1 = b$ is as before,
and we also know that the coefficient $b_2$ is only determined 
by $r$ and by the moments of $X$ up to order 6. In fact, one may evaluate 
$b_2$ on the basis of equality (7.5) of Proposition 7.2,
which specializes Proposition 5.1 to the range 
$6 \leq s < 8$. Since the formula for the coefficient $A_2 = A_2(r)$
is somewhat complicated, we will not go into tedious computations.

Now, from (1.5) it follows that
\bee
h_r(Z_{n+1}) - h_r(Z_n) 
 & = &
B_1\,\Big(\frac{1}{n} - \frac{1}{n+1}\Big) + 
b_2\,\Big(\frac{1}{(n+1)^2} - \frac{1}{n^2}\Big) + o(n^{-2}) \\
 & = &
\frac{B_1}{n(n+1)} + o(n^{-2}),
\ene
which thus proves Theorem 1.3 in case of finite $r$.

\vskip5mm
\section{{\bf Maximum of density (the case $r = \infty$)}}
\setcounter{equation}{0}

\vskip2mm
\noindent
Recall that $N_{\infty}(X) = \|p\|_{\infty}^{-2}$, when a random variable $X$
has density $p$. An expansion similar to the one of Proposition 5.1 can also 
be obtained for $\|p_n\|_{\infty}$ and hence for $N_\infty(Z_n)$. 
In order to deduce monotonicity, let us assume that $\beta_6 < \infty$.
 
From the non-uniform local limit theorem it follows that
$ \|p_n - \varphi_6\|_{\infty} = o(n^{-2})$ as $n \to \infty$,
where $\varphi_6$ is the Edgeworth expansion of order $6$. Hence
\be
\|p_n\|_{\infty} = \|\varphi_6\|_{\infty} + o(n^{-2}). 
\en
Here
$$ 
\varphi_6(x) = \varphi(x) 
\Big(1 + Q_1(x)\frac{1}{\sqrt{n}} + Q_2(x)\frac{1}{n} + 
Q_3(x)\frac{1}{n^{\frac{3}{2}}} + Q_4(x) \frac{1}{n^2} \Big), 
$$
where the polynomials $Q_k(x)$ are the same as in Section 3.

Let us find an asymptotic expansion for $\|\varphi_6\|_{\infty}$. 
Since $\varphi_6(x)$ is vanishing at infinity, there exists a point $x_6(n)$ 
such that $\|\varphi_6\|_{\infty} = |\varphi_6(x_6(n))|$. Since also
the functions $\varphi(x)\, Q_k(x)$ are bounded, we have
$|\varphi_6(x)| = O(\frac{1}{\sqrt{n}})$ uniformly in the region 
$|x| \geq \sqrt{\log n}$.
On the other hand, 
$$ 
\vp_6(0) = \vp(0) + \vp(0) \sum_{k=1}^4 Q_k(0)\,n^{-\frac{k}{2}} \geq 
\frac{1}{2}\,\vp(0)
$$
for $n$ large.
Therefore, $\vp_6(0) > |\vp_6(x)|$ for all $n$ large enough, as long as
$|x| \geq \sqrt{\log n}$, and we conclude that
\begin{eqnarray}
\|\vp_6 \|_{\infty} \, = \sup_{|x| \leq \sqrt{\log n}} |\vp_6(x)| \qquad
{\rm and} \qquad |x_6(n)| \leq \sqrt{\log n}.
\end{eqnarray}

Since $x=x_6(n)$ is the point of local extremum, we have $\vp_6'(x) = 0$, 
that is,
\be
x = 
\frac{Q_1'(x) - x Q_1(x)}{\sqrt{n}} + 
\frac{Q_2'(x) - x Q_2(x)}{n} + 
\frac{Q_3'(x) - x Q_3(x)}{n^{\frac{3}{2}}} + 
\frac{Q_4'(x) - x Q_4(x)}{n^2}. 
\en
Using (10.2), we deduce from $(10.3)$ that
$x_6(n) = O\big( \frac{1}{\sqrt{n}}\,(\log n)^{\frac{13}{2}}\big)$ and hence
$|x_6(n)| \leq 1$ for all $n$ large enough. But then, from $(10.3)$ again, 
$x_6(n) = O(\frac{1}{\sqrt{n}})$. For $x=x_6(n)$, we thus have
$$
\frac{x Q_3(x)}{n^{\frac{3}{2}}} = O\big(n^{-5/2}\big), \qquad
\frac{Q_4'(x)}{n^2} = O\big(n^{-5/2}\big), \qquad
\frac{x Q_4(x)}{n^2} = O\big(n^{-5/2}\big),
$$
and $(10.3)$ is simplified to 
$$
x = 
\frac{Q_1'(x) - x Q_1(x)}{\sqrt{n}} + 
\frac{Q_2'(x) - x Q_2(x)}{n} + 
\frac{Q_3'(x)}{n^{\frac{3}{2}}} + O\big(n^{-5/2}\big).
$$

The Chebyshev-Hermite polynomials satisfy the relation
$H_k'(x) - x H_k(x) = -H_{k+1}(x)$, so
\bee
H_3'(x) - x H_3(x) 
 & = &
- H_4(x) \ = \ -3 + 6x^2 - x^4 \\
H_4'(x) - x H_4(x)
 & = &
- H_5(x) \ = \ -15 x + 10 x^3 - x^5 \\
H_6'(x) - x H_6(x) 
 & = &
- H_7(x) \ = \ 105\,x - 105\, x^3 + 21\,x^5 - x^7.
\ene
Once $x=O(\frac{1}{\sqrt{n}})$, then
\bee
\frac{Q_1'(x) - x Q_1(x)}{\sqrt{n}} 
 & = &
\frac{\gamma_3}{3!}\, \frac{H_3'(x) - x H_3(x)}{\sqrt{n}} \\ 
 & = &
\frac{\gamma_3}{3!}\, \frac{-3 + 6x^2 - x^4}{\sqrt{n}} \ = \
- \frac{\gamma_3}{2\sqrt{n}} + \gamma_3 \frac{x^2}{\sqrt{n}} + 
O\big(n^{-5/2}\big)
\ene
and
\bee
\frac{Q_2'(x) - x Q_2(x)}{n} 
 & = &
\frac{\gamma_3^2}{2!\,3!^2}\, \frac{H_6'(x) - x H_6(x)}{n} + 
\frac{\gamma_4}{4!}\, \frac{H_4'(x) - x H_4(x)}{n} \\ 
 & = &
\frac{\gamma_3^2}{2!\,3!^2}\, \frac{105\,x - 105\, x^3 + 21\,x^5 - x^7}{n} + 
\frac{\gamma_4}{4!}\, \frac{-15 x + 10 x^3 - x^5}{n} \\
 & = &
\Big(\frac{105}{2!\,3!^2}\, \gamma_3^2 - \frac{15}{4!}\,\gamma_4\Big)\, 
\frac{x}{n} + O\big(n^{-5/2}\big).
\ene
Since
\bee
Q_3'(x) 
 & = &
\frac{1}{3!^4}\, \gamma_3^3\, H_9'(x) + 
\frac{1}{3!\,4!}\,\gamma_3 \gamma_4\, H_7'(x) + 
\frac{1}{5!}\,\gamma_5\, H_5'(x) \\
 & = &
\frac{9}{3!^4}\, \gamma_3^3\, H_8(x) + 
\frac{7}{3!\,4!}\,\gamma_3 \gamma_4\, H_6(x) + 
\frac{5}{5!}\,\gamma_5\, H_4(x) \\
 & = &
\frac{9 \cdot 105}{3!^4}\, \gamma_3^3 - 
\frac{7 \cdot 15}{3!\,4!}\,\gamma_3 \gamma_4 + 
\frac{5 \cdot 3}{5!}\,\gamma_5 + O(x^2),
\ene
we also have 
$$
\frac{Q_3'(x)}{n^{\frac{3}{2}}} = 
\Big(\frac{945}{3!^4}\, \gamma_3^3 - \frac{105}{3!\,4!}\,\gamma_3 \gamma_4 + 
\frac{15}{5!}\,\gamma_5\Big)\,\frac{1}{n^{\frac{3}{2}}} + 
O\big(n^{-5/2}\big).
$$
As a result,
\begin{eqnarray}
x \ = \ x_6(n) 
 & = & 
- \frac{\gamma_3}{2 \sqrt{n}} + \gamma_3 \frac{x^2}{\sqrt{n}} + 
\Big(\frac{105}{2 \cdot 3!^2}\, \gamma_3^2 - \frac{15}{4!}\, \gamma_4 \Big) 
\frac{x}{n} \nonumber \\ 
 & ~ & \qquad + \ 
\Big(\frac{945}{3!^4}\, \gamma_3^3 - \frac{105}{3! 4!}\, \gamma_3 \gamma_4 + 
\frac{15}{5!}\, \gamma_5\Big) 
\frac{1}{n^{\frac{3}{2}}} + O\big(n^{-5/2}\big).
\end{eqnarray}
One may use this asymptotic equation to find an expansion for $x_6(n)$
in powers of $1/\sqrt{n}$. Indeed, first we immediately obtain that
\begin{eqnarray*}
x = x_6(n) = 
- \frac{\gamma_3}{2 \sqrt{n}} + O\big(n^{-\frac{3}{2}}\big),
\end{eqnarray*}
implying
$$ 
\frac{x^2}{\sqrt{n}} = \frac{\gamma_3^2}{4} \frac{1}{n^{\frac{3}{2}}} + 
O\big(n^{-5/2}\big), \qquad
\frac{x}{n} = -\frac{\gamma_3}{2} \frac{1}{n^{\frac{3}{2}}} + 
O\big(n^{-5/2}\big).
$$
Inserting the above to (10.4), we deduce that
\begin{eqnarray*}
x = x_6(n) 
 & = & 
- \frac{\gamma_3}{2 \sqrt{n}} \\ 
 & & + \
\left(\frac{\gamma_3^3}{4} - \frac{\gamma_3}{2} 
\Big(\frac{105}{2 \cdot 3!^2}\, \gamma_3^2 - \frac{15}{4!}\, \gamma_4\Big) + 
\frac{945}{3!^4}\, \gamma_3^3 - \frac{105}{3! \cdot 4!}\, \gamma_3 \gamma_4 + 
\frac{15}{5!}\, \gamma_5 \right) \frac{1}{n^{\frac{3}{2}}} \\
 & & + \ 
O\big(n^{-5/2}\big) \\ 
 & = & 
\frac{a_1}{\sqrt{n}} + \frac{a_2}{n^{\frac{3}{2}}} + O\big(n^{-5/2}\big),
\end{eqnarray*}
with coefficients
$$
a_1 \, = \, -\frac{1}{2}\,\gamma_3, \qquad
a_2 \, = \, \frac{1}{4}\,\gamma_3^3 - \frac{5}{12}\, \gamma_3 \gamma_4 + 
\frac{1}{8}\,\gamma_5.
$$
In particular, $a_1 = a_2 = 0$ and therefore $x = x_6(n) = O\big(n^{-5/2}\big)$, 
as long as the distribution of $X$ is symmetric about the origin 
(in which case $\gamma_3 = \gamma_5 = 0$).

Still in the general case, keeping these coefficients, we deduce for
$x = x_6(n)$ that
\bee
x 
& = &
\frac{1}{\sqrt{n}}\,\Big(a_1 + a_2\,\frac{1}{n}\Big) + O\big(n^{-5/2}\big), 
\qquad
x^2 
\ = \
\frac{1}{n}\,\Big(a_1^2 + 2 a_1 a_2\,\frac{1}{n}\Big) + O\big(n^{-5/2}\big), \\
x^3
& = &
\frac{1}{n^\frac{3}{2}}\,a_1^3 + O\big(n^{-5/2}\big), \qquad
x^4
\ = \
\frac{1}{n^2}\,a_1^4 + O\big(n^{-5/2}\big), \qquad
x^p
 \, = \,
O\big(n^{-5/2}\big) \qquad (p \geq 5).
\ene
Hence
\bee
\frac{Q_1(x)}{\sqrt{n}} 
 & = &
\frac{\gamma_3}{6\sqrt{n}}\,(x^3 - 3x) \\
 & = &
\frac{\gamma_3}{6\sqrt{n}}\,\Big(\frac{1}{n^\frac{3}{2}}\,a_1^3 - 
\frac{3}{\sqrt{n}}\,\Big(a_1 + a_2\,\frac{1}{n}\Big)\Big) + 
O\big(n^{-5/2}\big)
 \ = \
\frac{\gamma_3^2}{4n} + \frac{b_1}{n^2} + O\big(n^{-5/2}\big)
\ene
with
$$ 
b_1 = \frac{\gamma_3}{3!}\, (a_1^3 - 3a_2). 
$$
Similarly,
\bee 
\frac{Q_2(x)}{n} 
 & = &
\Big(\frac{\gamma_3^2}{2!\,3!^2}\, H_6(x) + \frac{\gamma_4}{4!}\, H_4(x)\Big)\,
\frac{1}{n} \\
 & = &
\Big(\frac{\gamma_3^2}{2!\,3!^2}\, (-15 + 45x^2) + 
\frac{\gamma_4}{4!}\, (3 - 6x^2)\Big)\, \frac{1}{n} + O\big(n^{-5/2}\big)\\
 & = &
\Big(\frac{3}{4!}\, \gamma_4 - \frac{15}{2!\cdot 3!^2}\,\gamma_3^2\Big)\,\frac{1}{n} + 
\Big(\frac{45}{2!\cdot 3!^2}\,\gamma_3^2 - \frac{6}{4!}\, \gamma_4\Big)\,
\frac{x_2}{n} + O\big(n^{-5/2}\big) \\ 
 & = &
\Big(\frac{3}{4!}\, \gamma_4 - \frac{15}{2!\cdot 3!^2}\,\gamma_3^2\Big)\,\frac{1}{n} + 
\frac{b_2}{n^2} + O\big(n^{-5/2}\big)
\ene
with
$$ 
b_2 = 
\Big(\frac{45}{2 \cdot 3!^2}\, \gamma_3^2 - \frac{6}{4!}\,\gamma_4\Big)\, a_1^2. 
$$
Next,
\bee
\frac{Q_3(x)}{n^{\frac{3}{2}}} 
 & = &
\Big(\frac{\gamma_3^3}{3!^4}\, H_9(x) + 
\frac{\gamma_3 \gamma_4}{3!\,4!}\, H_7(x) + 
\frac{\gamma_5}{5!}\, H_5(x)\Big)\,\frac{1}{n^{\frac{3}{2}}} \\
 & = &
\Big(\frac{\gamma_3^3}{3!^4}\, 945\,x - 
\frac{\gamma_3 \gamma_4}{3!\,4!}\, 105\,x + 
\frac{\gamma_5}{5!}\, 15\,x\Big)\,\frac{1}{n^{\frac{3}{2}}} + 
O\big(n^{-5/2}\big)
 \ = \
\frac{b_3}{n^2} + O\big(n^{-5/2}\big)
\ene
with
$$ 
b_3 = \Big(\frac{945}{3!^4}\, \gamma_3^3 - 
\frac{105}{3! 4!}\, \gamma_3 \gamma_4 + \frac{15}{5!}\, \gamma_5\Big)\,a_1, 
$$
and finally
\bee
\frac{Q_4(x)}{n^2} 
 & = &
\Big(\frac{\gamma_3^4}{4!\,3!^4}\, H_{12}(x) + 
\frac{\gamma_3^2 \gamma_4}{2!\, 3!^2\,4!}\, H_{10}(x) + 
\frac{\gamma_3 \gamma_5}{3!\,5!}\, H_8(x) \\
 & & + \ 
\frac{\gamma_4^2}{2!\, 4!^2}\, H_8(x) + 
\frac{\gamma_6}{6!}\, H_6(x)\Big)\, \frac{1}{n^2} \\ 
 & = &
\Big(\frac{\gamma_3^4}{4!\,3!^4}\,10\,395 - 
\frac{\gamma_3^2 \gamma_4}{2!\, 3!^2\,4!}\,945 + 
\frac{\gamma_3 \gamma_5}{3!\,5!}\,105 \\
 & & + \ 
\frac{\gamma_4^2}{2!\, 4!^2}\,105 -
\frac{\gamma_6}{6!}\,15\Big)\, \frac{1}{n^2} + O\big(n^{-5/2}\big)
 \ = \
\frac{b_4}{n^2} + O\big(n^{-5/2}\big) \\
\ene
with
$$ 
b_4 = \frac{10\,395}{4!\cdot 3!^4} \gamma_3^4 - 
\frac{945}{2 \cdot 3!^2 4!}\, \gamma_3^2 \gamma_4 + 
\frac{105}{3!\cdot 5!}\, \gamma_3 \gamma_5 + 
\frac{105}{2 \cdot 4!^2}\,\gamma_4^2  -
\frac{15}{6!}\, \gamma_6. 
$$

Note that in the case of symmetric distributions,
$b_1 = b_2 = b_3 = 0$, while 
$
b_4 = \frac{105}{2 \cdot 4!^2}\,\gamma_4^2 - \frac{15}{6!}\, \gamma_6.
$

Now, as $x \to 0$,
$$
\frac{\vp(x)}{\|\vp\|_{\infty}} = 
1 - \frac{1}{2}\, x^2 + \frac{1}{8}\, x^4 + O(x^6), 
$$
and recall that, for $x = x_6(n)$, we have
$x^2 = \frac{1}{n}\,(a_1^2 + 2 a_1 a_2\,\frac{1}{n}) + O(n^{-5/2})$ and 
$x^4 = \frac{1}{n^2}\, a_1^4 + O(n^{-5/2})$. Thus,
$$ 
\frac{\vp(x)}{\|\vp\|_{\infty}} = 1 - \frac{a_1^2}{2n} + 
\Big(\frac{a_1^4}{8} - a_1 a_2\Big) \frac{1}{n^2} + O\big(n^{-5/2}\big). 
$$
Therefore, denoting $b = b_1 + b_2 + b_3 + b_4$, we get
\begin{eqnarray*}
\frac{\|\vp_6 \|_{\infty}}{\|\vp\|_{\infty}}
 & = & 
\frac{\vp_6(x)}{\|\vp\|_{\infty}} \\ 
 & = & 
\frac{\varphi(x)}{\|\vp\|_{\infty}}
\Big(1 + \frac{Q_1(x)}{\sqrt{n}} + \frac{Q_2(x)}{n} + 
\frac{Q_3(x)}{n^{\frac{3}{2}}} + \frac{Q_4(x)}{n^2}\Big) \\ 
 & = & 
\Big(1 - \frac{a_1^2}{2n} + \left( \frac{a_1^4}{8} - a_1 a_2 \right) 
\frac{1}{n^2}\Big)
\Big(1 + 
\Big(\frac{1}{4}\,\gamma_3^2 + \frac{3}{4!}\, \gamma_4 - 
\frac{15}{2!\cdot 3!^2}\,\gamma_3^2\Big)\,\frac{1}{n} + \frac{b}{n^2}\Big) + 
O\big(n^{-5/2}\big) \\ 
 & = &
1 + \Big(-\frac{1}{2}\, a_1^2 + \frac{1}{4}\,\gamma_3^2 + 
\frac{3}{4!}\, \gamma_4 - \frac{15}{2! \cdot 3!^2}\, \gamma_3^2\Big) 
\frac{1}{n} \\ 
 & & + \ 
\left(b + \frac{1}{8}\, a_1^4 - a_1 a_2 -
\frac{1}{2}\, \Big(\frac{1}{4}\,\gamma_3^2 + \frac{3}{4!}\, \gamma_4 - 
\frac{15}{2! \cdot 3!^2}\, \gamma_3^2\Big)\,a_1^2 \right) 
\frac{1}{n^2} + O\big(n^{-5/2}\big).
\end{eqnarray*}
Simplifying the term in front of $1/n$, we arrive at 
$$
\| \vp_6 \|_{\infty} = \| \vp\|_{\infty} + \frac{\|\vp\|_{\infty}}{n}\, A + 
\frac{\| \vp\|_{\infty}}{n^2}\,B + O\big(n^{-5/2}\big),
$$
where
\be
A = \frac{1}{8} \Big(\gamma_4 - \frac{2}{3}\, \gamma_3^2\Big), \qquad
B = b + \frac{1}{8}\,a_1^4 - a_1 a_2 -
\frac{1}{2}\, \Big(\frac{1}{4}\,\gamma_3^2 + \frac{3}{4!}\, \gamma_4 - 
\frac{15}{2! \cdot 3!^2}\, \gamma_3^2\Big)\,a_1^2.
\en

Using our assumptions, let us summarize by recalling the assertion 
(10.1): we get
\be
\|p_n\|_{\infty} = \|\vp\|_{\infty} \Big(1 + \frac{1}{n}\, A + 
\frac{1}{n^2}\,B\Big) + o\big(n^{-2}\big),
\en
where $A$ and $B$ are as above with
$$
a_1 \, = \, -\frac{1}{2}\,\gamma_3, \qquad
a_2 \, = \, \frac{1}{4}\,\gamma_3^3 - \frac{5}{12}\, \gamma_3 \gamma_4 + 
\frac{1}{8}\,\gamma_5.
$$
One can now reformulate this result in terms of the R\'enyi entropy
of index $r = \infty$. Since 
$N_\infty(Z_n) = \|p_n\|_\infty^{-2}$ and $N_\infty(Z) = \|\vp\|_\infty^{-2}$ 
for $Z \sim N(0,1)$, the expansion (10.6) yields:

\vskip5mm
{\bf Proposition 10.1.} {\sl If $\beta_6$ is finite, then as $n \to \infty$,
\be 
N_\infty(Z_n) \, = \,
N_\infty(Z) \Big(1 - \frac{\widetilde{A}}{n} + \frac{\widetilde{B}}{n^2}\Big) 
+ o\Big(\frac{1}{n^2}\Big)
\en
with
$$ 
\widetilde{A} = 
\frac{1}{4}\, \Big(\gamma_4 - \frac{2}{3}\, \gamma_3^2\Big), 
\qquad \widetilde{B} = 3A^2 - 2B,
$$
where the constants $A$ and $B$ are given in $(10.5)$.
}

\vskip5mm
{\bf Proof of Theorem 1.3 in case $r = \infty$.}
Denoting $\Delta_n = N_{\infty}(Z) -  N_\infty(Z_n)$, from (10.7) we get
$\Delta_{n+1} - \Delta_n = -\frac{\widetilde{A}}{n(n+1)} + o( \frac{1}{n^2})$.
\qed

\vskip5mm
In the case $\gamma_3 = \gamma_5 = 0$, for example when $X$ is symmetric, 
the coefficients in Proposition 10.1 are simplified.
Indeed, recalling the formula for $b_4$ in such a case, we have
$$ 
A = \frac{1}{8}\,\gamma_4, \qquad
B = b_4  = \frac{105}{2 \cdot 4!^2}\,\gamma_4^2 - \frac{15}{6!}\, \gamma_6,
$$
and therefore,
$$ 
\widetilde{A} = \frac{1}{4}\,\gamma_4, \qquad 
\widetilde{B} = 3A^2 - 2B = 
\frac{1}{24}\,\gamma_6 - \frac{13}{96}\,\gamma_4^2.
$$
As a consequence, the eventual monotonicity of $N_\infty(Z_n)$ can be 
deduced based on the sign of $\gamma_4$. However, if also $\gamma_4 = 0$, 
we need to look at the sign of $\gamma_6$.



\vskip5mm
\begin{thebibliography}{BH3}

\small
\bibitem[A-B-B-N]{A-B-B-N}
Artstein, S.; Ball, K. M.; Barthe, F.; Naor, A. Solution of Shannon's problem on the 
        monotonicity of entropy. J. Amer. Math. Soc.  17  (2004),  no. 4, 975--982. 

\vskip2mm 
\bibitem[B]{B} 
Barron, A. R. Entropy and the central limit theorem. Ann. Probab. 14 (1986), 
        no. 1, 336--342. 

\vskip2mm 
\bibitem[B-RR]{B-RR}
Bhattacharya, R. N.; Ranga Rao, R. Normal approximation and asymptotic expansions. 
        John Wiley \& Sons, Inc. 1976. Also: Soc. for Industrial and 
				Appl. Math., Philadelphia, 2010.

\vskip2mm
\bibitem[B-C]{B-C} 
Bobkov, S. G.; Chistyakov, G. P. Entropy power inequality for the R\'enyi entropy. 
        IEEE Trans. Inform. Theory  61 (2015),  no. 2, 708–-714.


\vskip2mm       
\bibitem[B-C-G1]{B-C-G1}
Bobkov, S. G.; Chistyakov, G. P.; G\"otze, F. Non-uniform bounds in local limit 
        theorems in case of fractional moments. I. Math. Methods Statist. 
		  	20 (2011), no. 3, 171--191; 20 (2011), no. 4, 269--287.

\vskip2mm       
\bibitem[B-C-G2]{B-C-G2}
Bobkov, S. G.; Chistyakov, G. P.; G\"otze, F. Rate of convergence and Edgeworth-type 
        expansion in the entropic central limit theorem. Ann. Probab. 41 (2013), no. 4, 
        2479--2512. 

\vskip2mm
\bibitem[B-M]{B-M} 
Bobkov, S. G.; Marsiglietti, A. Variants of the entropy power inequality. 
        IEEE Trans. Inform. Theory 63 (2017), no. 12, 7747--7752.

\vskip2mm
\bibitem[C-T]{C-T} 
Cover, T. M.; Thomas, J. A. Elements of information theory. Second edition. 
        Wiley-Interscience [John Wiley \& Sons], Hoboken, NJ, 2006, xxiv+748 pp.

\vskip2mm
\bibitem[D-C-T]{D-C-T} 
Dembo, A.; Cover, T. M.; Thomas, J. A. Information-theoretic inequalities. 
        IEEE Trans. Inform. Theory  37  (1991),  no. 6, 1501–-1518.
        
\vskip2mm
\bibitem[G-K]{G-K}         
Gnedenko, B. V.; Kolmogorov, A. N. Limit distributions for sums of independent random 
        variables. Translated and annotated by K. L. Chung. 
        With an Appendix by J. L. Doob. Addison-Wesley Publishing Company, Inc., 
        Cambridge, Mass., 1954. ix+264 pp. 
        
\vskip2mm
\bibitem[M-B]{M-B}         
Madiman, M.; Barron, A. Generalized entropy power inequalities and monotonicity 
        properties of information. IEEE Trans. Inform. Theory  53  (2007),  no. 7, 
        2317--2329.        

\vskip2mm
\bibitem[O-P]{O-P}
Osipov, L. V.; Petrov, V. V. On the estimation of the remainder term in the
        central limit theorem. Teor. Veroyatn. Primen. 12 (1967), 322--329.

\vskip2mm
\bibitem[P1]{P1}
Petrov, V. V. Local limit theorems for sums of independent random variables. 
        (Russian) Teor. Verojatnost. Primenen. 9 (1964), 343--352. 

\vskip2mm
\bibitem[P2]{P2}
Petrov, V. V. Sums of independent random variables. Translated from the Russian 
        by A. A. Brown. Ergebnisse der Mathematik und ihrer Grenzgebiete, 
			 Band 82. Springer-Verlag, New York-Heidelberg, 1975. x+346 pp. 
			 Russian ed.: Moscow, Nauka, 1972, 414 pp.			
					
					
\end{thebibliography}
\end{document}